\documentclass{amsart} 
\usepackage{amssymb, amscd, epsfig, supertabular, times}
\def\AA{{\mathbb A}} 
\def\BB{{\mathbb B}} 
\def\CC{{\mathbb C}} 
\def\DD{{\mathbb D}}

\def\HH{{\mathbb H}}

\def\LL{{\mathbb L}} 
 
\def\PP{{\mathbb P}} 
\def\QQ{{\mathbb Q}} 
\def\RR{{\mathbb R}}

\def\ZZ{{\mathbb Z}}

\def\IH{I\! H}

\def\G{{\Gamma}} 
\def\g{{\gamma}}

\def\eps{\varepsilon}

\def\st{{\rm st}} 
\def\sst{{\rm sst}}

\def\hyp{{\rm hyp}}

\def\bs{\backslash} 
 
\def\ii{{\sqrt{-1}}}

\def\Gcal{{\mathcal G}}
\def\Hcal{{\mathcal H}}

\def\Lcal{{\mathcal L}}
 
\def\Mcal{{\mathcal M}}
\def\Ocal{{\mathcal O}}

\def\Qcal{{\mathcal Q}}

\def\Scal{{\mathcal S}}

\def\pbr{{\rm PBr}}

\def\la{\langle} 
\def\ra{\rangle}
\def\half{\textstyle{\frac{1}{2}}}

\newcommand\bl{\operatorname{Bl}}
\newcommand\codim{\operatorname{codim}}
\newcommand\Div{\operatorname{Div}}

\newcommand\End{\operatorname{End}}

\newcommand\Hom{\operatorname{Hom}}

\newcommand\im{\operatorname{Im}}

\newcommand\Ker{\operatorname{Ker}} 
\newcommand\Coker{\operatorname{Coker}}

\newcommand\re{\operatorname{Re}}
\newcommand\res{\operatorname{Res}}

\newcommand\GL{\operatorname{GL}}
\newcommand\PGL{\operatorname{PGL}}
\newcommand\SL{\operatorname{SL}}

\newcommand\supp{\operatorname{supp}}

\newcommand\PU{\operatorname{PU}}

 \newtheorem{theorem}{Theorem}[section]
 \newtheorem{corollary}[theorem]{Corollary}
 \newtheorem{lemma}[theorem]{Lemma}
 \newtheorem{proposition}[theorem]{Proposition}
 \theoremstyle{definition}
 \newtheorem{definition}[theorem]{Definition}
 \theoremstyle{remark}
 \newtheorem{remark}[theorem]{Remark}
 \newtheorem{remarks}[theorem]{Remarks}
 \newtheorem*{example}{Example}
 \numberwithin{equation}{section}

\begin{document}

\title[Uniformization by Lauricella functions]
 {Uniformization by Lauricella functions---an overview of the theory of Deligne-Mostow}
\author[Eduard Looijenga]{Eduard Looijenga}

\address{Betafaculteit Universiteit Utrecht--Departement Wiskunde\\ 
Postbus 80.010\\ 
NL-3508 TA Utrecht\\ 
Nederland}

\email{looijeng@math.uu.nl}


\begin{abstract}
This is a survey of the Deligne-Mostow theory of Lauricella functions, or 
what almost amounts to the same thing,  of the period map for cyclic coverings of the Riemann sphere.
\end{abstract}

\maketitle

\section*{Introduction}
These notes are about the theory of hypergeometric functions in several variables. 
The functions in question generalize the 
Gau\ss \ hypergeometric function and are obtained as integrals of a multivalued 
differential of the form
\[
\eta_z:=(z_0-\zeta )^{-\mu_0}\cdots (z_n-\zeta )^{-\mu_n}d\zeta,
\]
where $z_0,\dots, z_n$ are pairwise distinct complex numbers (and are allowed to vary) and the exponents $\mu_k$ are taken in the open unit interval $(0,1)$
(and are always kept fixed).  If $\gamma$ is a path connecting some $z_k$ 
with some  $z_l$ whose relative interior avoids the $z_k$'s and if a determination 
of the differential along that path is chosen, then $\eta_z$ can be integrated
along $\g$ (the integral will indeed converge). That integral will depend 
holomorphically on $z=(z_0,\dots ,z_n)$, for if we vary $z$ a little, then we
can let $\g$ and the determination of $\eta_z$  follow this variation in a continuous manner. The  (multivalued) function of $z$ thus obtained is the type of hypergeometric function that takes the stage here. We now briefly explain which are the aspects of particular interest  that will make an appearance in this piece.

One readily finds that it is better not to focus on one such integral, but to consider 
all of them simultaneously, or rather, to consider for every $z$ as above
(and  fixed exponents),  the space $L_z$ of  power series 
expansions in $n+1$ complex variables at  $z$ that are linear 
combinations of such integrals. It turns out that this vector space $L_z$ has dimension $n$ and that the `tautological' map-germ $(\CC^{n+1},z)\to L_z^*$ sends $z$ to an element $\not= 0$ and has the following regularity property: if  $\Mcal_{0,n+2}$ stands for the configuration space of $(n+1)$-tuples in $\CC$ modulo affine-linear equivalence (which is also the configuration space of $(n+2)$-tuples on the Riemann sphere  modulo projective-linear equivalence),
then this map-germ drops to a local isomorphism $(\Mcal_{0,n+2}, [z])\to \PP(L_z^*)$. By analytic continuation we have an identification of $L_z$ with $L_{z'}$ for nearby $z'$ and the multivalued nature of the 
hypergeometric functions is reflected by the fact that if we let $z$ 
traverse a loop in the space of pairwise distinct $(n+1)$-tuples
and if let the elements of $L_z$  follow that loop by analytic continuation, 
then there results a linear (monodromy) transformation of $L_z$  which need not be the identity. The transformations of $L^*_z$ thus obtained  form a subgroup $\G$
of $\GL (L^*_z)$,  called the monodromy group of the system. The main questions addressed here are: 

1. When does $\G$ leave invariant a  Hermitian form
which is a positive definite, semidefinite or of  hyperbolic signature? 

2. In the situation of Question 1, when is $\G$ discrete as a subgroup of 
$\GL (L^*_z)$? (This is essentially equivalent to: when acts $\G$ 
properly on $\DD$?) And when is $\G$ arithmetic (in a naturally defined 
$\QQ$-algebraic group that contains $\G$)?

The answer to the first question  is short enough to give here: when 
$\mu_0+\cdots +\mu_n$ is $<1$, $=1$ or in the interval $(1,2)$ respectively (we are not claiming the converse).  In that case  $\Mcal_{0,n+2}$ acquires a metric of constant holomorphic curvature as follows. First, we observe that $\PP(L^*_z)$ contains  a complex symmetric manifold  of constant holomorphic curvature $\DD$ as an open subset on which $\G$ acts by isometries: we get respectively
all of $\PP(L^*_z)$ with its Fubini-Study metric, an affine space in 
$\PP(L^*_z)$ with a translation invariant metric or an open ball with its 
complex hyperbolic metric. But we also find that the local isomorphism $(\Mcal_{0,n+2},[z])\to \PP(L^*_z)$ lands in $\DD$, so that $\Mcal_{0,n+2}$ 
 inherits a metric from $\DD$.

Question 2 is harder. If $\G$ is discrete as well, then the exponents 
$\mu_k$ must be rational numbers. One of the main results 
states that  $\Mcal_{0,n+2}$ has then finite invariant 
volume and that its natural metric completion is an algebraic 
variety (we get a projective space in the elliptic and 
parabolic cases and in the hyperbolic case it is obtained by adding 
the stable orbits in a setting of geometric invariant theory). 
Deligne and Mostow gave sufficient conditions for discreteness, which 
were later weakened by Mostow and Sauter to make them sufficient as well.

If the $\mu_k$'s are rational, then there is the connection with the theory 
of period maps (regardless whether $\G$ is discrete): if $m$ is their smallest 
common denominator and if we write $\mu_k=d_k/m$, then the hypergeometric 
functions  become periods of the cyclic cover of $\CC$ defined by 
$w^m=(z_0-\zeta )^{d_0}\cdots (z_n-\zeta )^{d_n}$. For 
$\eta_z$ then  lifts to a regular univalued differential on this affine 
curve (regular resp.\ with simple poles at infinity when $\sum_k\mu_k$ 
is greater than resp.\ equal to $1$) and $\g$ is covered by a cycle such 
that the hypergeometric integral is the period of the lift over this cycle. 

As the reader will have gathered, this is mostly an account of work of 
Mostow (and his student Sauter) and 
of  Deligne-Mostow. It is self-contained in that the sense 
that we have included proofs  (except for a technical lemma needed 
for an arithmeticity criterion).  Occasionally our treatment somewhat 
differs from theirs. For instance, our discussion of invariant 
Hermitian forms does not use
the approach in \cite{delmost1} inspired by  Hodge theory, but 
rather follows 
the more pedestrian path in \cite{couw}. We also found it natural to use the 
language of orbifolds throughout. For some of the  history of the material 
expounded here, we  refer to the first and the last section  of 
\cite{delmost1} as well to the review \cite{ch}. In Section 
\ref{sect:other} we---very sketchily---mention some recent developments. 

This paper is based on a series of talks I gave at the CIMPA summer school (2005) 
in Istanbul. I thank my hosts, in particular Professor Uludag, for their hospitality and for making this summer school such a pleasant and fruitful experience.
 
\tableofcontents

\section{The Lauricella differential}\label{sect:lauricella}

\subsection{Definition and first properties}\label{subsect:first}
Assume given  real numbers $\mu_0,\dots ,\mu_n$ in the interval $(0,1)$, where $n>0$. We shall refer to the $(n+1)$-tuple $\mu=(\mu_0,\dots ,\mu_n)$ as a \emph{weight system} and we call its sum  $|\mu|:=\sum_{i=0}^n \mu_i$ 
the \emph{total weight} of $\mu$.  The \emph{Lauricella differential} of weight $\mu$ is
\[
\eta_z:=(z_0-\zeta )^{-\mu_0}\cdots (z_n-\zeta )^{-\mu_n}d\zeta,
\quad \text{with  } z=(z_0,\dots ,z_n)\in (\CC^{n+1})^\circ.
\]
(Here  $(\CC^{n+1})^\circ$ stands for the set of  $(z_0,\dots ,z_n)\in\CC^{n+1}$ whose components are pairwise  distinct.)
Athough this differential is multivalued, it has a natural determination
on a left half space  by taking there the 
value of the integrand whose argument is $<\pi/|\mu|$ in absolute value.
We further note that $\eta_z$ is locally integrable as a multivalued function: near
$z_k$, $\eta_z$ is  of the form $(\zeta-z_k)^{-\mu_k}\exp(\textit{holom})d\zeta$;
this is  the differential of a function of the form 
$\textit{const}+(\zeta-z_k)^{1-\mu_k}\exp(\textit{holom})$ and since  $1-\mu_k>0$, that function
takes a well-defined value in $z_k$. This implies that $\eta_z$ can be 
integrated along every 
\emph{relative arc} of $(\CC ,\{ z_0,\dots ,z_n\})$; by the latter we mean an 
oriented piecewise differentiable arc in $\CC$ whose end points lie in 
$\{ z_0,\dots ,z_n\}$, but which does not meet this set elsewhere. 

The behavior of the differential  at infinity is studied by means of the substitution $\zeta=\omega^{-1}$; this gives
\[
\eta_z=-(\omega z_0-1)^{-\mu_0}\cdots (\omega z_n-1)^{-\mu_n}\omega^{|\mu |-2}d\omega,
\]
which suggests to put $z_{n+1}:=\infty$ and $\mu_{n+1}:=2-|\mu |$. In case
$\mu_{n+1}<1$ (equivalently, $|\mu|>1$), $\eta_z$ is also (multivalued) integrable at $z_{n+1}$.

\begin{remark}
Following Thurston \cite{thurstontable}, we may think of $\eta_z$ as a way of putting a flat Euclidean structure on $\PP^1$ with singularities at $z_0,\dots ,z_{n+1}$: a local integral of
$\eta_z$ defines a metric chart with values in $\CC$, but now regarded as the Euclidean plane (so the associated metric is simply $|\eta_z|^2$). At $z_k$, $k\le n$, the metric space is isometric to a Euclidean cone
with total angle $2\pi (1-\mu_k)$; this is also true for $k=n+1$ in case 
$\mu_{n+1}<1$, or equivalently, $|\mu|>1$;
if $|\mu|=1$ resp.\ $|\mu|<1$, then a punctured  neighborhood of $\infty$ is isometric to a flat cylinder resp.\ the complement of a compact subset of a Euclidean 
cone with total angle $1-|\mu|$.
\end{remark}
 
Let be given relative arc $\gamma_z$ of $(\CC ,\{ z_0,\dots ,z_n\})$ and 
let also be given a determination of $\eta_z$ on $\gamma_z$ 
so that $\int_{\gamma_z} \eta$ is defined. Choose an open disks $D_k$ about 
$z_k$ in $\CC$ such that the $D_0,\dots ,D_n$ are pairwise disjoint. 
Then we can find for every $z'\in 
D_0\times\cdots\times D_n$, a relative arc $\gamma_{z'}$ of 
$(\CC ,\{ z'_0,\dots ,z'_n\})$  and a determination of $\eta_{z'}$ on 
$\supp (\gamma_{z'})$ such that both  depend continuously on $z'$ and 
yield the prescribed value for $z=z'$. Any primitive of $\eta$ near 
$(z,z_k)$ with respect to its second variable is (as a function of $(z',\zeta)$) of the form $g(z')+(\zeta-z'_k)^{1-\mu_k}h(\zeta,z')$, 
with $g$ and $h$ holomorphic and so the function
\[
z'\in D_0\times\cdots\times D_n\mapsto \int_{\gamma_{z'}} \eta_{z'}\in\CC
\]
is holomorphic. We call such a function (or
some analytic extension of it) a \emph{Lauricella function}. The Lauricella  functions (with given weight system $\mu$) span a complex vector space. We denote the space of germs of holomorphic functions at $z\in (\CC^{n+1})^\circ$ that are  germs of Lauricella 
functions by $L_z$. It is clear that for $z'\in D_0\times\cdots \times D_n$, we can naturally identify $L_{z'}$ with $L_z$. 
Here are some elementary properties of Lauricella functions (the proofs are left to the reader, who should be duely  careful with exchanging differention and integration in the proof of (c) ).

\begin{proposition}\label{prop:elementary} Any $f\in L_z$
\begin{enumerate}
\item[(a)] is translation invariant: $f(z_0+a,\dots ,z_n+a)=f(z_0,\dots ,z_n)$ for small $a\in\CC$,
\item[(b)] is homogeneous of degree $1-|\mu |$: $f(e^tz_0,\dots , e^tz_n)=e^{(1-|\mu |)t}f(z_0,\dots ,z_n)$ for small $t\in\CC$  and 
\item[(c)] obeys the system of differential equations
\[
\frac{\partial^2 f}{\partial z_k\partial z_l}=\frac{1}{z_k-z_l}\left( \mu_l
\frac{\partial}{\partial z_k}-\mu_k\frac{\partial f}{\partial z_l}\right),\quad 0\le k<l\le n.
\]
\end{enumerate}
\end{proposition}

The translation invariance of the Lauricella functions suggests to introduce
\[
V_n:= \CC^{n+1}/\text{main diagonal}\quad\text{and}\quad 
V_n^\circ:=(\CC^{n+1})^\circ/\text{main diagonal},
\]
as they are in fact defined on $V_n^\circ$. The homogeneity implies that when $|\mu |=1$, these functions are also constant on the 
$\CC^\times$-orbits and hence factor through $\PP (V_n^\circ)$; for reasons
which will become clear later, we call this the \emph{parabolic} case.

An important consequence of part (c) of the preceding proposition is 

\begin{corollary}\label{cor:inj}
The map which assigns to $f\in L_z$ its $1$-jet at $z$ is injective.
\end{corollary}
\begin{proof}
If $f\in L_z$, then its partial derivatives $f_k:=\frac{\partial f}{\partial z_k}$ satisfy
the system of ordinary differential equations
\[
\frac{\partial f_k}{\partial z_l}=\frac{1}{z_k-z_l}\left( \mu_l
f_k-\mu_kf_l\right),\quad k\not= l.
\]
We can complete this system as to get also such  equations for 
$\frac{\partial f_k}{\partial z_k}$ by using the fact $\sum_k f_k=0$ (which follows
from the translation invariance). The elementary theory of such systems says
that there is precisely one solution for it, once the initial conditions
$f_k(z)$ are prescribed. To such a solution corresponds  at most one element 
of $L_z$ up to a constant.
\end{proof}

\subsection{Lauricella arc systems}\label{subsect:arc}
\begin{definition}
Given  $(z_0,\dots ,z_n)\in\CC^{n+1}$, we define an \emph{L-arc system} 
as to be an oriented arc in the Riemann sphere $\PP^1=\CC\cup\{\infty\}$ from $z_0$ to $z_{n+1}=\infty$  which passes successively through $z_1,\dots z_n$ and
follows near $\infty$  the real axis in the positive direction. 
If $\delta$ is such an $L$-arc system, then we denote 
the  piece connecting $z_{k-1}$ with $z_k$ by $\delta_k$ and we often 
let $\delta$ also stand for the system of arcs $(\delta_1,\dots ,\delta_{n+1})$. 
\end{definition} 

The complement of the support of $\delta$ is simply connected and so we have a
well-defined determination of $\eta_z$ on this complement which extends the
one we already have on a left half space. We also extend $\eta_z$ to the 
support of $\delta$ itself by insisting that $\eta_z$ be continuous 
`from the left' (which makes the determination of
$\eta_z$ discontinuous along $\delta$). With these conventions,
$\int_{\delta_k}\eta_z$ has for $k=1,\dots ,n$ a well defined meaning (and
also makes sense for $k=n+1$ in case $\mu_{n+1}<1$). If we let $z$ vary in a small neighborhood, we get
an element of $L_z$ that we simply denote by $\int_{\delta_k}\eta $.
We denote by $\delta_k^-$ the arc connecting $z_{k-1}$ with $z_k$ that is 
`infinitesimally' to the right of $\delta_k$. By this we really mean that 
$\eta_z$  is given on  $\delta_k^-$ the determination  it gets as a limit from the right.
Notice that $\eta_z| \delta_k^-= \exp (-2\pi \ii (\mu_0+\cdots +\mu_{k-1}))\eta_z| \delta_k$.

\begin{theorem}\label{theorem:basis}
The functions $\int_{\delta_1}\eta,\dots, \int_{\delta_n}\eta$ define a basis  for 
$L_z$. Moreover, $L_z$ contains the constant functions if and only if we are in the parabolic case: $|\mu |=1$.
\end{theorem}
\begin{proof}
Any relative arc of $(\CC,\{z_0,\dots ,z_n\})$ is homotopic to a 
composite of the arcs $\delta_k$, and their inverses (we want 
the homotopy be such that the determination of $\eta$ varies continuously).
Since any two determinations of $\eta$ differ by a constant factor, 
this implies that the functions $\int_{\delta_1}\eta,\dots, \int_{\delta_n}\eta$ generate $L_z$.

If $|\mu|=1$, then $\eta_z$ is near $\infty$ equal to $-\zeta^{-1}d\zeta$. 
So then for a loop $\gamma$ which encircles $z_0,\dots ,z_n$ in the clockwise direction, we have
\[
\int_\gamma \eta_z=\int_\gamma -\zeta^{-1}d\zeta=2\pi \ii,
\]
which proves that $L_z$ contains the constant $2\pi\ii$.

It remains to show that if $a_1,\dots ,a_k,c\in\CC$ are such that
$\sum_{k=1}^n a_k \int_{\delta_k}\eta =c$, then $c\not=0$ implies $|\mu |=1$ and
$c=0$ implies that all $a_i$ vanish as well. We prove this with induction on $n$.
To this end, we consider a curve $z(s)$ in $(\CC^{n+1})^\circ$ of the form
$(z_0,\dots ,z_{n-2},0, s)$, with $s>0$ and an $L$-arc system $\delta (s)$
for $z(s)$ with $\delta_1,\dots ,\delta_{n-1}$ fixed and $\delta_n=[0, s]$.
By analytic continuation we may assume that 
$\sum_{k=1}^{n-1} a_k \int_{\delta_k}\eta_{z(s)} +a_n \int_{0}^{s}\eta_{z(s)} =c$. We multiply this identity with $s^{\mu_n}$ and investigate
what happens for $s\to\infty$. For $k<n$, 
\[
s^{\mu_n}\int_{\delta_k}\eta_{z(s)}=
\int_{\delta_k} (z_0-\zeta)^{-\mu_0}\cdots (z_{n-2}-\zeta)^{-\mu_{n-2}}(-\zeta)^{-\mu_{n-1}}(1-s^{-1}\zeta)^{-\mu_n}d\zeta,
\]
which for $s\to\infty$ tends to $\int_{\delta_k}\eta_{z'}$, where $z'=(z_0,\dots ,z_{n-1})$.
On the other hand, 
\begin{multline*}
\int_0^s(z_0-\zeta)^{-\mu_0}\cdots (-\zeta)^{-\mu_{n-1}}
(s-\zeta)^{-\mu_n}d\zeta\\
=s (-s)^{-|\mu|}\int_0^1 (-s^{-1}z_0+\zeta)^{-\mu_0}\cdots (\zeta)^{-\mu_{n-1}}
(-1+\zeta)^{-\mu_n}d\zeta\\
=s(-s)^{-|\mu|} +o(|s|^{1-|\mu|}), \quad  s\to\infty.
\end{multline*}
So we find that 
\[
s^{\mu_n}\left(c+a_n\left((-s)^{1-|\mu|}+o(|s|^{1-|\mu|}\right)\right)=\sum_{k=1}^{n-1} a_k\int_{\delta_k}\eta_{z'},
\quad s\to\infty .
\]
This shows that $c\not=0$ implies $|\mu |=1$ (and $a_n=(-1)^{-|\mu|}$). 
Suppose now $c=0$. If $\mu_n<|\mu|-1$, then the left hand side tends to zero
as $s\to\infty$ and so the right hand side must be zero. Our induction hypothesis
then implies that $a_1=\cdots =a_{n-1}=0$ and from this we see  that $a_n=0$, too.
If $\mu_n>|\mu|-1$, then we clearly must have $a_n=0$ and the induction hypothesis
implies that $a_1=\cdots =a_{n-1}=0$, also.
\end{proof}

\begin{remark}
So the space of solutions of the system of differential equations
in  Proposition \ref{prop:elementary}-c is in the nonparabolic case equal to $L_z\oplus\CC$, and contains $L_z$ as a hyperplane in the parabolic case.
\end{remark}

\subsection{The rank of the Schwarz map}\label{subsect:schwarz}
We find it convenient to modify our basis of Lauricella functions
by a scalar factor by putting
\begin{align*}
F_k(z,\delta):=&\int_{\delta_k} (\zeta-z_0)^{-\mu_0}\cdots  (\zeta-z_{k-1})^{-\mu_{k-1}}
 (z_k-\zeta)^{-\mu_k}\cdots  (z_n-\zeta)^{-\mu_n}d\zeta\\
=&\bar w_k\int_{\delta_k} \eta_z,\quad \text{where } w_k:=e^{\ii\pi (\mu_0+\cdots +\mu_{k-1})}.
\end{align*}
The notation now also displays the fact the value of the integral depends on the whole $L$-arc system
(which was needed to make $\eta_z$ univalued) and not just on $\delta_k$. Notice that if $z=x$ is real and $x_0<x_1<\cdots <x_n$ and $\delta$ consists of real intervals, then the integrand is real valued and positive and hence so is $F_k$.
Let us also observe that
\[
\int_{\delta_k}\eta_z=w_kF_k(z,\delta)\; \text{ and }\; \int_{\delta^-_k}\eta_z=\bar w_kF_k(z,\delta),
\]
where  the second identity follows from  the fact that 
$\eta_z|\delta^-_k=\bar w_k^2\eta_z|\delta_k$. 
So if we are in the parabolic case,
then the integral of $\eta_z$ along a clockwise  loop which encloses 
$\{ z_0,\dots ,z_n\}$ yields the identity $\sum_{k=1}^n (w_k-\bar w_k)F_k(z,\delta)=2\pi \ii$,
or equivalently,
\begin{equation}\label{eq:par}
\sum_{k=1}^n \im (w_k)F_k(z,\delta)=\pi.
\end{equation}
In other words, $F=(F_1,\dots ,F_n)$ then maps to the affine hyperplane
$\AA^{n-1}$ in $\CC^n$ defined by this equation.

\begin{corollary}\label{cor:lociso}
If  we are not in the parabolic case, then  $F=(F_1,\dots ,F_n)$, 
viewed as a multivalued map from
$V_n^\circ$  to $\CC^n$, is a local isomorphism and never 
takes the origin of $\CC^n$ as value. In the parabolic case,  
$F=(F_1,\dots ,F_n)$ factors through a local (multivalued) isomorphism 
from $\PP(V_n^\circ)$ to the affine hyperplane $\AA^{n-1}$ in $\CC^n$ 
defined by $\sum_{k=1}^n \im (w_k)F_k=\pi$. 
\end{corollary}
\begin{proof}
Given $(z,\delta)$, consider the $n$ covectors
$dF_1(z,\delta),\dots ,dF_n(z,\delta)$ in the cotangent space of $z$. 
According to corollary \ref{cor:inj}, a linear relation among them must arise
from a linear relation among the function germs $F_1,\dots ,F_n\in L_z$ and 
the constant function $1$. According to Theorem \ref{theorem:basis}, such  
a relation exists if and only if $|\mu|=1$. The corollary easily follows, 
except perhaps the claim that
$F$ is nowhere zero. But if $F_k(z,\delta)=0$ for all $k$, then we must have 
$|\mu |\not= 1$; since $F$ will be constant zero on  the $\CC^\times$-orbit 
through $z$  this contradicts the fact that $F$ is a local isomorphism.
\end{proof}

\begin{definition}
We call the multivalued map $F$ from  $V_n^\circ$ to $\CC^n$ the 
\emph{Lauricella map} and its projectivization $\PP F$ from  $\PP(V_n^\circ)$ to $\PP^{n-1}$ the \emph{Schwarz map} for the weight system $\mu$. 
\end{definition}

The above corollary tells us that the Schwarz map always is a local isomorphism (which in the parabolic case takes values in the affine open $\AA^{n-1}\subset\PP^{n-1}$).

\subsection{When points coalesce}\label{subsect:coalesce}
We investigate the limiting behavior of $F$ when some of the $z_k$'s come together. 
To be specific,  fix $0<r<n$ and let for $0<\eps<1$, $z_\eps=(\eps z_0,\dots \eps z_r,z_{r+1},\dots z_n)$ and see what happens when  $\eps\to 0$. 
We assume here that $z_1,\dots ,z_r$  lie inside the unit disk, whereas the others are outside that disk and choose $\delta$  accordingly: 
$\delta_1,\dots ,\delta_r$  resp.\ $\delta_{r+2},\dots ,\delta_{n+1}$ lie inside resp.\ outside the unit disk. 

Put $\mu':=(\mu_0,\dots ,\mu_r)$,
$z'=(z_0,\dots ,z_r)$. Then
\begin{align*}
F_k(z_\eps ,\delta)=\bar w_k\int_{\delta_k} (\varepsilon z_0-\zeta )^{-\mu_0}\cdots
(\varepsilon z_r-\zeta )^{-\mu_r} (z_{r+1}-\zeta )^{-\mu_{r+1}}\cdots
(z_n-\zeta )^{-\mu_n} d\zeta\\
=\varepsilon^{1-|\mu'|} 
\bar w_k\int_{\delta_k} (z_0-\zeta )^{-\mu_0}\cdots
(z_r-\zeta )^{-\mu_r} (z_{r+1}-\varepsilon\zeta )^{-\mu_{r+1}}\cdots
(z_n-\varepsilon\zeta )^{-\mu_n} d\zeta,
\end{align*}
where  in the last line (involving the passage to $\eps\zeta$ as  the new integration varable) $\delta_k$ must be suitably re-interpreted. So  for $k\le r$,
\begin{equation}\label{eq:normal}
\varepsilon^{|\mu'|-1}F_k(z_\eps,\delta)=
(1+ O(\varepsilon))
z_{r+1}^{-\mu_{r+1}}\cdots z_n^{-\mu_n} F'_k(z',\delta'),
\end{equation}
where $F'_k$ is a component of the Lauricella map with weight system $\mu'$:
\[
 F'_k(z',\delta')=\bar w_k \int_{\delta_k} (z_0-\zeta )^{-\mu_0}\cdots (z_k-\zeta )^{-\mu_k}d\zeta.
 \]  
If  $k>r$ and  in case $k=r+1$, $|\mu'|<1$, we find
\begin{equation}\label{eq:longi}
F_k(z_\eps,\delta)=(1+ O(\varepsilon))\bar w_k\int_{\delta_k} (-\zeta )^{-|\mu'|}
(z_{r+1}-\zeta )^{-\mu_{r+1}}\cdots (z_n-\zeta )^{-\mu_n} d\zeta.
\end{equation}
Assume now $|\mu'|<1$. Then these estimates suggest to replace in $F=(F_1,\dots ,F_n)$, for $k\le r$ , $F_k$ by  $\varepsilon^{|\mu'|-1}F_k(z,\delta)$. In geometric terms, this amounts to 
enlarging the domain and range of $F$: now view it as a multivalued map defined an open subset of the blowup $\bl_{(z_0,\dots ,z_r)}V_n$ of the diagonal defined $z_0=\cdots =z_r$  and as mapping to  the blowup $\bl_{(F_1,\dots ,F_r)}\CC^n$ of the subspace 
of $\CC^n$ defined by $F_1=\cdots =F_r=0$. It maps the 
exceptional divisor (defined by $\varepsilon=0$)
to the exceptional divisor $\PP^{r-1}\times \CC^{n-r}\subset \bl_{(F_1,\dots ,F_r)}\CC^n$. If we identify the exceptional divisor in the domain
with $\PP(V_r)\times V_{1+n-r}$  (the second component begins with the  common value of $z_0,\dots,z_r$), then we see that the first component of this restriction is the Schwarz map $\PP F'$ for the weight system $\mu'$ and the second component is 
$\bar w_{r}$ times the Lauricella map for the weight system $(\mu',\mu_{r+1},\dots ,\mu_n)$. 

If several such clusters are forming, then we have essentially a product situation.

We shall also need to understand what happens  when $|\mu'|=1$. Then taking the limit for $\eps\to 0$ presents a problem for $F_{r+1}$ only (the other components have well-defined limits). This is related to the fact that $\eta_z$ is univalued on the unit circle $S^1$; by the theory of  residues we then have
\[
\int_{S^1} \eta_z= \int_{S^1} 
(z_0-\zeta )^{-\mu_0}\cdots (z_n-\zeta )^{-\mu_n} d\zeta=
2\pi\ii z_{r+1}^{-\mu_{r+1}}\cdots z_n^{-\mu_n}.
\]
We therefore replace $\eta_z$ by $\hat\eta_z:=z_{r+1}^{\mu_{r+1}}\cdots z_n^{\mu_n}\eta_z$ and $F$  by  $\hat F:=z_{r+1}^{\mu_{r+1}}\cdots z_n^{\mu_n}F$. This does not change the Schwarz map, of course. Notice however, that  now $\int_{S^1} \hat\eta_z=2\pi\ii$.

\begin{lemma}\label{lemma:parlimit}
Assume that $\mu'$ is of parabolic type: $|\mu'|=1$.
Define Lauricella data $\mu'':=(\mu_{r+1},\dots ,\mu_{n+1})$, $z'':=(z_{r+1}^{-1},\dots ,z_n^{-1},0)$ and let $\delta''=(\delta''_1,\dots ,\delta''_{n-r})$ be  the image of $(\delta_{r+2},\dots ,\delta_{n+1})$ under the map $z\mapsto z^{-1}$.  Then we have
\[
\hat F_k(z_\eps,\delta)=
\begin{cases}
(1+ O(\varepsilon))F'_k(z',\delta') &\text{ when $1\le k\le r$,}\\
(1+ O(\varepsilon))F''_{k-r-1}(z'',\delta'') &\text{ when $r+2\le k\le n$,}
\end{cases}
\]
whereas $\lim_{\eps\to 0}\re \hat F_{r+1}(z_\eps,\delta)=+\infty$. 
Moreover,  $\sum_{k=1}^r\im (w_k)\hat F_k(z,\delta)=\pi$.
\end{lemma}
\begin{proof}
The assertion for $k\le r$ is immediate from our previous calculation. 
For $1\le i\le n-r-1$  we find
\begin{multline*}
\hat F_{r+1+i}(z_\eps,\delta)=\\
=w_{r+1+i}\int_{\delta_{r+1+i}} (\eps z_0-\zeta )^{-\mu_0}\cdots
(\varepsilon z_r-\zeta )^{-\mu_r} (1-\frac{\zeta}{z_{r+1}} )^{-\mu_{r+1}}\cdots
(1-\frac{\zeta}{z_n} )^{-\mu_n} d\zeta\\
=-w''_i \int_{\delta''_i} (\eps z_0-\zeta^{-1} )^{-\mu_0}\cdots
(\eps z_r-\zeta^{-1} )^{-\mu_r} (1-\frac{1}{\zeta z_{r+1}})^{-\mu_{r+1}}\cdots
(1-\frac{1}{\zeta z_n} )^{-\mu_n}\frac{d\zeta}{-\zeta^2}\\
=w''_i \int_{\delta''_i} (1-\eps z_0\zeta)^{-\mu_0}\cdots
(1-\eps z_r\zeta )^{-\mu_r} (\frac{1}{z_{r+1}}-\zeta)^{-\mu_{r+1}}\cdots
(\frac{1}{z_n}-\zeta )^{-\mu_n}(-\zeta)^{-\mu_{n+1}}d\zeta\\
=(1+ O(\varepsilon))F''_i(z'',\delta'').
\end{multline*}
As to the limiting behavior of $\hat F_{r+1}$, observe that
\begin{multline*}
\hat F_{r+1}(z_\eps,\delta)=\\
=-\int_{\delta_{r+1}} (\eps z_0-\zeta )^{-\mu_0}\cdots
(\varepsilon z_r-\zeta )^{-\mu_r} (1-\frac{\zeta}{z_{r+1}} )^{-\mu_{r+1}}\cdots
(1-\frac{\zeta}{z_n} )^{-\mu_n} d\zeta\\
=\int_{\delta_{r+1}} (\zeta-\varepsilon z_0)^{-\mu_0}\cdots
(\zeta-\varepsilon z_r)^{-\mu_r} (1-\frac{\zeta}{z_{r+1}} )^{-\mu_{r+1}}\cdots
(1-\frac{\zeta}{z_n} )^{-\mu_n} d\zeta.
\end{multline*}
For $\eps\to 0$, the integrand tends to
 $\zeta^{-1} (1-\zeta/z_{r+1} )^{-\mu_{r+1}}\cdots(1-\zeta/z_n )^{-\mu_n}$, from
which the asserted limiting behaviour easily follows.
The last assertion follows from the fact that  $\int_{S^1} \hat\eta_z=2\pi\ii$
(see the derivation of Equation \eqref{eq:par}).
 \end{proof}

So if we regard the Schwarz map as defined on an open subset of 
$\bl_{(z_0,\dots ,z_k)}\PP (V_n)$, then its composite with 
the projection
of $\PP^n\to \PP^{n-1}$ obtained by omitting $F_r$ is on the exceptional 
divisor  given by $[F'_1:\cdots :F'_r:F''_1:\cdots :F''_{n-1-r}]$.

\subsection{Monodromy group and monodromy cover}\label{subsect:monodromy}
We begin with making a few remarks about the fundamental group of 
$(\CC^{n+1})^\circ$. We take $[n]=(0,1,2,\dots ,n)$ as a base point for 
$(\CC^{n+1})^\circ$ and use the same symbol for its image in $V_n^\circ$. 
The projection $(\CC^{n+1})^\circ\to V_n^\circ$ induces an 
isomorphism on fundamental groups:   
$\pi_1((\CC^{n+1})^\circ,[n])\cong  \pi_1(V_n^\circ, [n])$. 
This group is known as the \emph{pure} (also called \emph{colored}) 
\emph{braid group with $n+1$ strands}; we denote it by
$\pbr_{n+1}$. Another characterization of $\pbr_{n+1}$ is that
as the group of connected components of the group of diffeomorphisms 
$\CC\to\CC$ that are the identity outside a compact subset of $\CC$ and 
fix each $z_k$. 

If $\alpha$ is a path in $(\CC^{n+1})^\circ$ from $z$ to $z'$, and if we are
given an $L$-arc system $\delta$ for  $z$, then we can carry that system
continuously along when we traverse $\alpha$; we end up with an 
$L$-arc system $\delta'$ for  $z'$ and this $L$-arc system  will 
be unique up to isotopy. In this way $\pbr_{n+1}$  acts on the set 
of isotopy classes of  $L$-arc systems. 
It is not hard to see that this action is principal: 
for every ordered pair of isotopy classes of $L$-arc systems, 
there is a unique element of $\pbr_{n+1}$ which carries the 
first one onto the second one. 

The group $\pbr_{n+1}$ has a set of distinguished elements, called 
\emph{Dehn twists},
defined as follows.  The basic Dehn twist is a diffeomorphism of the 
annulus $D_{1 ,2}\subset\CC: 1\le |z|\le 2$; it is defined by 
$re^{\ii \theta}\mapsto re^{\ii (\theta +\phi (r))}$, where $\phi$ 
is a differentiable function which is zero resp.\ $2\pi$ on a 
neighborhood of $1$ resp.\ $2$ (all such diffeomorphisms of $D_{1,2}$ 
are isotopic relative to the boundary $\partial D_{1,2}$). 
If $S$ is an oriented surface, 
and we are given an orientation preserving diffeomorphism 
$h:D_{1,2}\to S$, then the Dehn twist on the
the image and the identity map on its complement define a diffeomorphism of $S$,
which is also called a Dehn twist. Its isotopy class only depends on the 
isotopy class of the image of the counter clockwise oriented unit circle 
(as an oriented submanifold of $S$). These embedded circles occur here as 
the isotopy classes of embedded circles in $\CC-\{z_1,\dots ,z_n\}$. 
A particular case of interest is 
such a circle encloses precisely two points of  $\{z_1,\dots ,z_n\}$, 
say $z_k$ and $z_l$.
The isotopy class of such a circle defines and is defined by the isotopy 
class of an unoriented path in $\CC-\{z_1,\dots ,z_n\}$ that connects $z_k$ 
and $z_l$ (the boundary of a regular neighborhood of such a path gives  
an embedded circle). The element  of the pure braid group associated to this 
is called \emph{simple}; if we choose for every  pair $0\le k<l\le n$ a 
simple element, then the resulting collection of simple elements is known 
to generate $\pbr_{n+1}$.

\smallskip
There is a standard way to obtain a covering of $V_n^\circ$ on which 
$F$ is defined as a univalued map. Let us recall this in the present case. 
First notice that if $\alpha$ is a path in $(\CC^{n+1})^\circ$ from $z$ to $z'$, 
then analytic continuation
along this path gives rise to an isomorphism of vector spaces 
$\rho_\mu(\alpha): L_z\to L_{z'}$. This is compatible with composition: 
if $\beta$ is a path in $(\CC^{n+1})^\circ$ from $z'$ to $z''$, then 
$\rho_\mu (\beta)\rho_\mu (\alpha)=\rho_\mu (\beta\alpha)$ (we use 
the functorial convention for composition of paths: $\beta\alpha$ means 
$\alpha$ followed by $\beta$). A loop in $(\CC^{n+1})^\circ$ based at $[n]$
defines an element $\rho_\mu (\alpha)\in \GL (L_{[n]})$ and we thus get a 
representation $\rho_\mu $ of $\pbr_{n+1}$ in $L_{[n]}$. The image of this  \emph{monodromy representation} is called the   
the \emph{monodromy group} (of the Lauricella system with weight system 
$\mu$); we shall denote that group by $\Gamma_\mu$, or simply by $\G$. 
The monodromy representation defines a 
$\G$-covering $\widetilde{V_n^\circ}$ of $V_n^\circ$ on which the $F_k$'s 
are univalued.
A point of $\widetilde{V_n^\circ}$ can be represented as a pair
$(z,\alpha)$, where $\alpha$ is a path in $\CC^{n+1}$ from $[n]$ to $z$, 
with the understanding that $(z',\alpha')$ represents the same point if 
and only if $z-z'$ lies on the main diagonal (so that $L_{z'}=L_z$) and 
$\rho_\mu(\alpha)=\rho_\mu (\alpha')$. The action
of $\G$ on $\widetilde{V_n^\circ}$ is then  given as follows: if 
$g\in \G$ is represented by the loop $\alpha_g$ in $\CC^{n+1}$ from $[n]$, 
then
$g.[(z,\alpha)]=[(z,\alpha \alpha_g^{-1})]$. But it is often more useful to
represent a point of $\widetilde{V_n^\circ}$  as a pair $(z,\delta)$, where 
$\delta$
is an $L$-arc system for $z$, with the understanding that $(z',\delta')$ represents the same point if and only if $z-z'$ lies on the main diagonal and $F_k(z,\delta)=F_k(z',\delta')$ for all $k=1,\dots ,n$. 
For this description we see right away that the basic Lauricella functions define a univalued holomorphic map 
\[
F=(F_1,\dots F_n): \widetilde{V_n^\circ}\to \CC^n.
\]
Since $[(z,\delta)]$ only depends on the isotopy class of $\delta$, the action of $\G$ is also easily explicated in terms of the last description. The germ of $F$ at the base point defines an isomorphism $L^*_{[n]}\cong \CC^n$: $c=(c_1,\dots ,c_n)\in\CC^n$ defines the linear form on $L_z$ which sends $F_k$ to $c_k$. If we let $\G$ act on
$\CC^n$ accordingly (i.e., as the dual of $L_{[n]}$), then $F$ becomes $\G$-equivariant.

The $\CC^\times$-action on $V_n^\circ$ given by scalar multiplication will lift
not necessarily to a $\CC^\times$-action on $\widetilde{V_n^\circ}$, but to one of a (possibly) infinite covering
$\widetilde{\CC^\times}$. For this action, $F$ is homogeneous of degree $1-|\mu|$.
Let us denote by  $\PP(\widetilde{V_n^\circ})$ the $\widetilde{\CC^\times}$-orbit space of $\widetilde{V_n^\circ}$.

\subsection{Invariant Hermitian forms}\label{subsect:hermform}
Our goal is to prove the following theorem.

\begin{theorem}\label{theorem:metric}
If $|\mu |<1$, then the monodromy group $\G$ leaves invariant a positive definite Hermitian form $H$ on $\CC^n$.

If $|\mu |=1$ (the parabolic case), then $\G$ leaves invariant  a positive definite Hermitian form $H$ on the (linear) translation hyperplane of the affine hyperplane 
$\AA^{n-1}$ in $\CC^n$, defined by $\sum_{k=1}^n \im (w_k)F_k=0$.

If $1<|\mu |<2$, then the monodromy group $\G$ leaves invariant a hyperbolic Hermitian form $H$ on $\CC^n$ (so of signature $(n-1,1)$) with the property that
$H(F(\tilde z,\tilde z))<0$ for all $\tilde z\in \widetilde{V_n^\circ}$.
\end{theorem}

Before we begin the proof, let us make the following observation.
If $W$ is a finite dimensional complex vector space, then by definition a point 
$u$ of $\PP(W)$ is given by a one-dimensional subspace $L_p\subset W$. 
An exercise shows that the complex tangent space $T_p\PP(W)$ of $\PP (W)$ at $p$ is naturally isomorphic to $\Hom (L_p,W/L_p)$. If we are also given a Hermitian form $H$ on $W$ which is nonzero on $L_p$, then it determines a Hermitian form $H_p$ on
$T_p\PP(W)\cong \Hom (L_p,W/L_p)$ as follows: since $H$ is nonzero on $L_p$, 
the $H$-orthogonal complement $L_p^\perp$ maps isomorphically $W/L_p$; 
if we choose a generator $u\in L_p$  and think of a tangent vector as a 
linear map $\phi :L_p\to L_p^\perp$, then we put $H_p(\phi ,\phi'):
=|H(u,u)|^{-1}H(\phi (u),\phi'(u))$. This
is clearly independent of the generator $u$. It is also clear that $H_p$ only
depends on the conformal equivalence class of $H$: it does not change if we
multiply $H$ by a positive scalar. 

If $H$ is positive definite, then so is 
$H_p$ for every $p\in \PP (W)$. In this way $\PP(W)$ acquires a Hermitian metric, 
known as the \emph{Fubini-Study} metric. It is in fact a K\"ahler manifold on which
the unitary group of $(W,H)$ acts transitively.

There is another case of interest, namely when $H$ has hyperbolic signature: if
we restrict ourselves to the set $\BB (W)$ of $p\in\PP(W)$ for which $H$ is 
negative on $L_p$, then $H_p$ is positive definite as well. This defines a metric on $\BB(W)$ which is 
invariant under the unitary group of $(W,H)$. If we choose a basis of linear forms
$u_0,\dots ,u_m$ on $W$ such that $H$ takes the standard form
$H(u,u)=-|u_0|^2+|u_1|^2+\cdots + |u_m|^2$, then we see that $\BB(W)$ is defined
in $\PP(W)$ by the inequality $|u_1/u_0|^2+\cdots +|u_m/u_0|^2<1$, which is simply the open unit ball in complex $m$-space. We call $\BB(W)$ a \emph{complex-hyperbolic space} and the metric defined above, the \emph{complex-hyperbolic metric}. As in the Fubini-Study case,
this metric makes $\BB(W)$ a  K\"ahler manifold on which the unitary group of $(W,H)$ acts transitively. For $m=1$ we recover the complex unit disk with its Poincar\'e metric.

Returning to the situation of Theorem \ref{theorem:metric}, we see that in all three cases
$\PP F$ is a local isomorphism mapping to a homogeneous  K\"ahler manifold:
when $|\mu|<1$, the range is a Fubini-Study space $\PP_{n-1}$ (this notatation is a private one: the  subscript is supposed to distinguish it from the  metricless projective space $\PP^{n-1}$), 
for $|\mu |=1$ we get a complex affine space
with a translation invariant metric (indeed, denoted here by $\AA_{n-1}$) and 
when $|\mu|>1$ we get  a complex ball  $\BB_{n-1}$ with its complex-hyperbolic metric. Since these structures are $\G$-invariant, we can state this more poignantly:  the weight system $\mu$ endows $\PP(V_n^\circ)$ with a natural K\"ahler metric locally
isometric with a Fubini-Study metric, a flat metric or a complex-hyperbolic metric.
We will therefore use  the corresponding terminology for the  cases $|\mu|<1$ and $1<|\mu|<2$  and call them the \emph{elliptic} and \emph{hyperbolic} case, respectively.

\medskip
Theorem \ref{theorem:metric} follows from more specific result that takes a bit of preparation to formulate. We shall associate to the weight system $\mu$ a  Hermitian form $H$ on $\CC^n$ or on the hyperplane in $\CC^n$ defined by $\sum_{k=1}^{n} \im(w_k)F_k=0$ in $\CC^n$, depending on whether $|\mu|$ is integral.
We do this somewhat indirectly. Let  $\tilde H$ be the Hermitian form on $\CC^{n+1}$ defined by 
\[
\tilde H(F,G)=\sum_{1\le j<k\le n+1}  \im(w_j\bar w_k)F_k\bar G_j.
\]
The $\tilde H$-orthogonal complement in $\CC^{n+1}$ of the last basis vector 
$e_{n+1}$ is the hyperplane $\tilde A$ defined by $\sum_{k=1}^{n+1} \im(w_k)F_k=0$. 
When $|\mu|\not\in\ZZ$, the projection $\tilde A\subset\CC^{n+1}\to\CC^n$ 
(which forgets the last coordinate) is an isomorphism  
(since  $w_{n+1}=e^{\pi\ii |\mu|}$, $\im (w_{n+1})\not=0$ in that case) and thus 
identifies $\CC^n$ with this $\tilde A$; we let $H$ then be the restriction 
of $\tilde H$ to $\CC^n$. If $|\mu|\in\ZZ$, then $\im (w_{n+1})=0$ and hence the  
projection $\tilde A\subset\CC^{n+1}\to\CC^n$ has kernel $\CC e_{n+1}$
and image the hyperplane $A$ in $\CC^n$ defined by $\sum_{k=1}^{n} \im(w_k)F_k=0$. 
So then $\tilde H$ induces a Hermitian form on $A$. The following proposition implies Theorem \ref{theorem:metric}.

\begin{proposition}
The form $H$ is $\G$-invariant for all weight systems $\mu$.
For $0<|\mu|\le 1$, the form $H$ is positive definite. For
$1<|\mu|<2$, $H$ is of hyperbolic signature and we have $H(F(z,\delta),F(z,\delta))=N(z)$, where 
\[
N(z)=-\frac{\ii}{2}\int_\CC \eta\wedge\bar\eta=-\int_\CC
|z_0-\zeta |^{-2\mu_0}\cdots |z_n-\zeta |^{-2\mu_n}| d\rm{(area)}.
\]
\end{proposition}
\begin{proof}
The assertions about the signature of $H$ involve a linear algebra calculation
that we leave to the reader (who may consult \cite{chl}).
We do the hyperbolic case first, so assume $1<|\mu |<2$. First notice that the integral
defining $N(z)$ converges (here we use that $|\mu|>1$) and takes on
a value which is real and negative. We claim that
\begin{equation}
N(z)=\sum_{1\le j<k\le n+1} w_j\bar w_k \bar F_j(z,\delta)F_k(z,\delta).
\end{equation}
To see this, let us integrate $\eta=\eta_z$, using the determination defined by $\delta$:
$\Phi (\zeta):=\int_{z_0}^\zeta \eta$, where the path of integration is not allowed
to cross $\supp (\delta)$. We have $d\Phi=\eta$ outside $\supp (\delta)$ and 
by Stokes theorem
\[
N(z)=-\frac{\ii}{2} \int_\CC \eta\wedge\bar\eta=
\frac{\ii}{2} \int_\CC d(\bar\Phi\eta) =
\frac{\ii}{2}\sum_{k=1}^{n+1}\left(
\int_{\delta_k} \bar\Phi\eta -\int_{\delta^-_k} \bar\Phi\eta\right).
\]
As to the last terms, we observe that  on $\delta_k$ we
have $\Phi (\zeta)=\sum_{j<k} w_jF_j+\int_{z_{k-1}}^\zeta \eta$
(we abbreviate $F_j(z,\delta)$ by $F_j$), where the 
last integral is taken over a subarc of $\delta_k$. Likewise, on $\delta^-_k$: 
$(\Phi |{\delta_k^-})(\zeta)=\sum_{j<k} \bar w_jF_j+\int_{z_{k-1}}^\zeta \bar w_k^2\eta$. Hence on $\delta_k$ we have 
\begin{multline*}
\bar\Phi\eta - (\bar\Phi\eta|{\delta_k^-})
=\sum_{j<k} \left(\bar w_j\bar F_j+\int_{z_{k-1}}^\zeta \bar\eta\right)\eta-
\sum_{j<k} \left(w_j\bar F_j+\int_{z_{k-1}}^\zeta w_k^2\bar \eta\right)\bar w_k^2\eta=\\
=\sum_{j<k} \left(\bar w_j -w_j\bar w_k^2\right)\bar F_j \eta
=\sum_{j<k} \left(\bar w_j w_k -w_j\bar w_k\right)\bar F_j \bar w_k\eta,
\end{multline*}
which after integration over $\delta_k$ yields
\[
\int_{\delta_k} \bar\Phi\eta -\int_{\delta^-_k} \bar\Phi\eta=
\sum_{j<k} \left(\bar w_jw_k -w_j\bar w_k\right)\bar F_j F_k
=\frac{2}{\ii}\sum_{j<k} \im(w_j\bar w_k)\bar F_j F_k.
\]
Our claim follows if we substitute this identity in the formula for $N$ above.

We continue the proof.  The claim implies that $H(F(z,\delta),F(z,\delta))=N(z)$.
The function $N$ is obviously $\G$-invariant (it does not involve $\delta$).
Since $N$ determines $H$, so is $H$. So this settles the hyperbolic case.

For the elliptic and parabolic cases we may verify by hand that it is invariant under a generating set of monodromy transformations, but a computation free argument, based analytic continuation as in \cite{chl}, is perhaps more satisfying. It runs as follows: if we choose a finite set of generators $\alpha_1,\dots ,\alpha_N$ of $\pbr_{n+1}$, then for every weight system $\mu$ we have a projective linear transformation 
$\PP\rho_\mu(\alpha_i)$ of $\PP^{n-1}$ that depends in a real-analytic manner on $\mu$. We will see that the Hermitian forms $h_\mu$ defined on an open subset of the tangent bundle of $\PP^{n-1}$ also depend real-analytically on $\mu$; so 
if $h_\mu$ is preserved by the $\PP\rho_\mu (\alpha_i)$'s for a nonempty open 
subset of $\mu$'s, then it is preserved for all weight systems for which this makes sense.
Hence $\PP\rho_\mu(\alpha_i)$ multiplies $H$  by a scalar. For $1<|\mu|<2$ this scalar is constant 1. Another analytic continuation argument implies that it is 1 for all $\mu$. 
\end{proof}

\subsection{Cohomological interpretation via local systems of rank one}\label{subsect:localsystem}
 
We sketch a setting in terms of which  the Hermitian form $H$ is best understood. It will not play a role in what follows (hence may be skipped), although it will reappear in a more conventional context (and formally  independent of this discussion) in Section \ref{sect:modular}. The reader should consult \S 2 of \cite{delmost1} for a more thorough treatment.

Fix complex numbers $\alpha_0,\dots,\alpha_{n}$ in $\CC^\times$.
Let $\LL$ be a local system of rank one
on $U:=\CC-\{ z_0,\dots ,z_n\}=\PP^1-\{ z_0,\dots, z_{n+1}\}$ such that the (counterclockwise) monodromy around $z_k$ is multiplication by $\alpha_k$.  It is unique up to isomorphism. We fix  a nonzero multivalued section $e$ of $\LL$ by choosing a nonzero section of $\LL$ on some left half plane and then extend that section to the universal cover of $U$ (defined by that left half plane). 
Denote by $\Lcal:=\Ocal_U\otimes_\CC\LL$ the underlying holomorphic line bundle. 
So if $\mu_k\in\CC$ is such that $\exp(2\pi\mu_k\ii)=\alpha_k$, then $s(\zeta):=\prod_{k=1}^n (z_k-\zeta)^{-\mu_k}\otimes e$ can be understood as a generating section of $\Lcal$. Likewise, $s d\zeta$ is a generating of $\Omega(\Lcal)=\Omega_U\otimes_\CC\LL$. Notice that $\Lcal$ comes with a connection $\nabla: \Lcal\to \Omega(\Lcal)$ characterized by
\[
\nabla (s)=\left(\sum_{k=0}^n \frac{\mu_k}{z_k-\zeta}\right)s d\zeta
\]
and that $\LL$ is recovered from the pair $(\Lcal, \nabla)$ as the kernel of $\nabla$.

The topological Euler characteristic of a rank one local system on a space homotopy equivalent to a finite cell complex is independent of that local system and hence
equal to the topological Euler characteristic of that space.
So the topological  Euler characteristic of $\LL$ is $-n$. Now assume that $\alpha_k\not=1$ for all $k$. This ensures that $\LL$ has no nonzero section. 
As there is no cohomology in degrees $\not=0,1$, this implies that
$\dim H^1(\LL)=n$. Moreover, if $j:U\subset \PP^1$ is the inclusion, then 
the stalk of $j_*\LL$ in $z_k$ is represented by the sections of $\LL$ on a punctured neighborhood of $z_k$, hence is zero unless $k=n+1$ and $\alpha_0\cdots \alpha_n=1$: then it is nonzero. So the map of complexes $j_!\LL\to j_*\LL$ has cokernel a one-dimensional skyscraper sheaf at $\infty$ or is an isomorphism. This implies that for the natural map $i: H^1_c(\LL)\to H^1(\LL)$,
$\dim\Ker (i)=\dim\Coker (i)$  is $1$ or $0$, depending on whether or not $\alpha_0\cdots \alpha_n=1$. We denote the image of $i$ by $\IH^1(\LL)$.

A relative arc $\alpha$ plus a section of $\LL^\vee$ over its relative  interior defines
a relative cycle of $(\PP^1,\{z_0,\dots ,z_{n+1}\})$ with values in $\LL^\vee$
and hence an element $[\alpha]$ of the relative homology space $H_1(\PP^1, \{z_0,\dots, z_{n+1}\};\LL^\vee)$. 
Alexander duality identifies the latter cohomology space with the dual of $H^1(\LL)$. 
To make the connection with the preceding, let us identify $\eta$ with $sd\zeta$
(we need not assume  here  that $\mu_k\in (0,1)$), so that we have a De Rham class $[\eta]\in H^1(\LL)$. 
If  we are given an $L$-arc system $\delta$ and choose the  determination of $e$ on $\delta_k$ prescribed by the arc system, then $\{\bar w_k[\delta_k]\}_{k=1}^n$ is a basis of $H_1(\PP^1, \{z_0,\dots, z_{n+1}\};\LL^\vee)$ and the value of 
$[\eta]$ on $\bar w_k[\delta_k]$ is just $F_k(z,\delta)$.

We have a perfect (Poincar\'e) duality 
$H^1_c(\LL)\times  H^1(\LL^\vee)\to \CC$, which, if  cohomology is represented by means of forms,  is given by integration over $U$ of the cup product.
Suppose now in addition that $|\alpha_k|=1$ for all $k$. This implies that $\LL$ carries a flat metric; we choose this metric to be the one for which $e$ has unit length.
The metric may be viewed as a $\CC$-linear isomorphism
of sheaves $\overline{\LL}\to \LL^\vee$ (here $\overline{\LL}$ stands for the local system $\LL$ with its complex conjugate complex structure) so that our perfect duality becomes a bilinear map $H^1_c(\LL)\times  \overline{H^1(\LL)}\to \CC$. 
We multiply that map by $\half\ii$ and denote the resulting sesquilinear map
$h:H^1_c(\LL)\times  H^1(\LL)\to \CC$.  Then $h$ is Hermitian in the sense that if $\alpha,\beta\in H^1_c(\LL)$, then $h(\alpha, i_*\beta)= \overline{h(\beta, i_*\alpha)}$, in particular, it induces a nondegenerate Hermitian form on $\IH ^1(\LL)$. This is just minus the form we defined in Subsection \ref{subsect:hermform}. If we take $\mu_k\in (0,1)$ for $k=0,\dots , n$ 
and assume $1<|\mu|<2$ (so that $\mu_{n+1}\in (0,1)$ also and $i$ is an isomorphism), then $h([\eta],[\eta])$ equals $\half\ii\int_{\CC}\eta\wedge \bar\eta$ indeed and hence equals $-N(z) = -H(F(z),F(z))$.

\section{Orbifolds and discrete monodromy groups}\label{sect:orbifolds}

\subsection{Monodromy defined by a simple Dehn twist}\label{subsect:dehn}
Let be given a relative arc $\gamma_0$ in $(\CC ,\{ z_0,\dots ,z_n\})$ which connects $z_k$ with
$z_l$, $k\not= l$. This defines a Dehn twist $D(\gamma_0)$ and hence an element \
$T$ of $\pbr_{n+1}$. We determine the action of $T$ on $\CC^n$.
For this we need to make $\eta_z$ univalued. Suppose we are given a straight piece of arc $\gamma_1$ that begins in $z_l$, but is otherwise disjoint from $\gamma_0$ so
that a neighborhood of $\supp(\gamma_0)$ minus $\supp(\gamma_0\gamma_1)$
is simply connected. Then choose a determination for $\eta_z$
on this simply connected open subset and let $\eta_z|\gamma_0$ be the limit from the left. Let  $\gamma$ resp.\ $\gamma'$ be an arc which ends in $z_k$ resp.\ $z_l$, but otherwise avoids  $\{ z_0,\dots ,z_n\}\cup\supp (\gamma_0\gamma_1)$
(we also assume that $\gamma'$ stays on the right of
$\gamma_1\gamma_0$).  Then from a picture one sees that
\begin{align*}
\int_{T(\gamma)}\eta_z&=\int_{\gamma}\eta_z + (1-w_l^2)\int_{\gamma_0} \eta_z,\\
\int_{T(\gamma)}\eta_z&=\int_{\gamma}\eta_z + (-w_l^2+w_k^2w_l^2)\int_{\gamma_0} \eta_z,\\
\int_{T(\gamma_0)}\eta_z&= w_k^2w_l^2\int_{\gamma_0} \eta_z.
\end{align*}
Remembering that $w_k^2w_l^2=e^{2\pi\ii (\mu_k+\mu_l)}$, one easily deduces from these formulae:

\begin{corollary}\label{cor:dehnmon}
If $\mu_k+\mu_l\not=1$, then $T$ acts in $\CC^n$ semisimply as a complex reflection over an angle $2\pi (\mu_k+\mu_l)$. If $\mu_k+\mu_l=1$, then $T$ acts in $\CC^n$ as a nontrivial unipotent transformation.
In particular, $T$ acts with finite order if and only if $\mu_k+\mu_l$ is a rational number
$\not= 1$.  
\end{corollary}

By a  \emph{complex reflection} we mean here a semisimple transformation which fixes a 
hyperplane pointwise. In the elliptic and hyperbolic cases, $T$ will be an orthogonal reflection with respect the Hermitian form $H$; in the parabolic case, it will be restrict to  $\AA_{n-1}$ as an orthogonal  affine reflection.

\subsection{Extension of the evaluation map}\label{subsect:extension}

The $\G$-covering $\widetilde{V_n^\circ}\to V_n^\circ$ can sometimes be extended
as a ramified $\G$-covering over a bigger open subset $V_n^f\supset V_n^\circ$ of $V_n$ (the superscript $f$ stands for $f$inite ramification; we may write $V_n^{f(\mu)}$ instead of $V_n^f$  if such precision is warranted). This means that we find a normal analytic variety $\widetilde{V_n^f}$ which contains $\widetilde{V_n^\circ}$ as an open-dense subset and to which the $\G$-action extends such that the $\G$-orbit space can be identified with $V_n^f$. This involves a standard tool
in analytic geometry that presumably goes back to Riemann and now falls under the heading of \emph{normalization}. It goes like this.
If $v\in V_n$ has a connected neighborhood $U_v$ in $V_n^\circ$
such that  one  (hence every) connected component of its preimage in  $\widetilde{V_n^\circ}$ is finite over $U_v\cap V_n^\circ$, then the $\G$-covering over
$U_v\cap V_n^\circ$ extends to a ramified $\G$-covering over $U_v$. The property
imposed on $U_v$ is equivalent to having finite  monodromy over $U_v\cap V_n^\circ$.
The extension is unique and so if $V_n^f$ denotes the set of $v\in V_n$ with this property, then a 
ramified $\G$-covering $\widetilde{V_n^f}\to V_n^f$ exists as asserted. 
The naturality of the construction also ensures that the $\tilde\CC^\times$-action on $\widetilde{V_n^\circ}$ (which covers the $\CC^\times$-action on $V_n^\circ$) extends to $\widetilde{V_n^f}$.

The space $V_n$ receives a natural stratification from the stratification of $\CC^{n+1}$ by its diagonals and since the topology of $V_n^f$ along strata does not change,
$V_n^f$ is an open union of strata. The codimension one strata are of the form
$D_{k,l}$, $0\le k<l\le n$, parameterizing the $z$ for which $z_k=z_l$, but no other equality among its components holds.

\begin{lemma}\label{lemma:ellhyperplane}
The stratum $D_{k,l}$ lies in $V^f_n$ if and only if $\mu_k+\mu_l$ is a rational number $\not=1$. The Schwarz map extends across  the preimage of  $\PP (D_{k,l})$) holomorphically if and only if $\mu_k+\mu_l<1$ and it does so as a local isomorphism
if and only if $1-\mu_k-\mu_l$ is the reciprocal of a positive integer. If $|\mu|\not=1$, then the corresponding assertions also hold for the Lauricella map. 
\end{lemma}
\begin{proof}
In order that $D_{k,l}\subset V_n^f$, it is necessary and sufficient that we have finite
monodromy along a simple loop around $D_{k,l}$. This monodromy is the image of a Dehn twist along a circle separating $z_k$ and $z_l$ from the other elements of
$\{z_0,\dots ,z_n\}$. So the first assertion follows from Corollary \ref{cor:dehnmon}. 

If $\gamma_0$ connects $z_k$ with $z_l$ within the circle specified above, then $\int_{\gamma_0}\eta_z=(z_k-z_l)^{1-\mu_k-\mu_l}\exp({\text{holom}})$. This  is essentially a consequence  of the identity
\[
\int_0^\varepsilon t^{-\mu_k}(t-\varepsilon)^{-\mu_l}dt=
\varepsilon^{1-\mu_k-\mu_l}\int_0^1 t^{-\mu_k}(t-1)^{-\mu_l}dt.
\]
Suppose now that $\mu_k+\mu_l\in\QQ-\{ 1\}$ and write $1-\mu_k-\mu_l=p/q$ with $p,q$ relatively prime integers with $q>0$. So the order of the monodromy is $q$
and over the preimage of a point of $D_{k,l}$, we have a coordinate $\tilde z_{k,l}$ with the property that $z_k-z_l$ pulls back to $\tilde z_{k,l}^q$. Hence $\int_{\gamma_0}\eta_z$ pulls back to $\tilde z_{k,l}^p$. In order that the Schwarz map extends over the preimage of $D_{k,l}$ holomorhically (resp.\ as a local isomorphism), a necessary condition is that the Lauricella function $\int_{\gamma_0}\eta_z$ (which after all may be taken as part of a basis of Lauricella functions) is holomorphic (resp.\ has a nonzero derivative everywhere). This means that $p>0$ (resp.\ $p=1$). It is not hard to verify that this is also sufficient.
\end{proof}

\subsection{The elliptic and parabolic cases}\label{subsect:ellpar}
Here the main result is:

\begin{theorem}[Elliptic case]\label{thm:ell}
Suppose that $|\mu |<1$ and that for all $0\le k<l\le n$,  $1-\mu_k-\mu_l$ is the reciprocal of an integer. Then $\G$ is a finite complex reflection group in $\GL(n,\CC)$ (so that in particular $V_n^f=V_n$) and $F: \widetilde{V_n}\to \CC^n$ is a $\G$-equivariant isomorphism which drops to an isomorphism $V_n\to \G\bs\CC^n$. 
\end{theorem}

So $\PP(V_n)$ acquires in these cases the structure of an
orbifold modeled on Fubini-Study space.
At the same time we prove a proposition that will be also useful later.
Observe that stratum of $V_n$ is given by a partition of $\{0,\dots ,n\}$: 
for $z$ in this stratum we have $z_k=z_l$ if and only if $k$ and $l$ 
belong to the same part. Let us say that this  stratum is 
\emph{stable relative to $\mu$} if its associated partition has the property 
that every part has $\mu$-weight $<1$. We denote by
$V_n^\st\subset V_n$  (or $V_n^{\st(\mu)}\subset V_n$) the union of 
stable strata.

\begin{proposition}\label{prop:elliptic}
Suppose that whenever $0\le k<l\le n$ are such that  $\mu_k+\mu_l<1$,   then $1-\mu_k-\mu_l$ is the reciprocal of an integer. Then $V_n^{\st}\subset V_n^f$, $\widetilde{V_n^{\st}}$ is a complex manifold. The Lauricella map extends holomorphically over this manifold and has the same regularity properties as  the map it extends: it is a local isomorphism when we are not in the parabolic case, whereas in the parabolic case, the Schwarz map defines a local isomorphism to $\AA_{n-1}$.  
\end{proposition}

We shall need: 

\begin{lemma}\label{lemma:covering}
Let $f:X\to Y$ be a local diffeomorphism from a manifold to a connected Riemannian manifold. Assume that $X$ is complete for the induced metric.
Then $f$ is a covering map.
\end{lemma}
\begin{proof}
We use the theorem of Hopf-Rinow which says that completeness is equivalent  to the property that every geodesic extends indefinitely as a geodesic.
Let $y\in Y$. Choose  $\eps>0$ such that the $\eps$-ball $B(y,\eps)$ is the diffeomorphic  image of the $\eps$-ball in $T_yY$ under  the exponential map. It is enough to show that every $x\in f^{-1}B(y,\eps)$ has a neighborhood which is mapped by $f$ diffeomorphically onto $B(y,\eps)$.
Since $X$ is complete, there is a (geodesic) lift of the geodesic in $B(y,\eps)$ from $f(x)$ to $y$ which begins in $x$. Then the end point $x_0$ of that lift lies in $f^{-1}y$. Then $B(x_0\eps)$  contains $x$ and maps diffeomorphically onto $B(y,\eps)$.
\end{proof}

We now begin the proofs  of Theorem \ref{thm:ell} and Proposition 
\ref{prop:elliptic}. Let us write  $A_k$ for the assertion of Theorem \ref{thm:ell}  for $k+1$ points and $B_k$ for the assertion of Proposition \ref{thm:ell} for 
elliptic strata of codimension $\le k$. Let us observe that $B_1$ holds: an elliptic stratum of codimension one is a stratum of the form $D_{k,l}$ satisfying the hypotheses of Lemma \ref{lemma:ellhyperplane}. We now continue with induction following the scheme below.

\begin{proof}[Proof that $A_k$  implies $B_k$]
Consider a stratum of codimension $k$. Let us first assume that it is irreducible in the sense that it is given by a single part. Without loss of generality we may then assume
that it is the open-dense in the locus $z_0=\cdots =z_k$. This is the situation we studied in Subsection \ref{subsect:coalesce} (mainly for this reason, as we can now confess). We found that $F$ extends to as a multivalued map defined on an open subset of the blowup $\bl_{(z_0,\dots ,z_k)}V_n$ going to the blowup $\bl_{(F_1,\dots ,F_k)}\CC^n$. On the the exceptional divisor,  $F$ is the product of the
Schwarz map for $\mu'=(\mu_0,\dots ,\mu_k)$ and the Lauricella map for
$(|\mu'|, \mu_{k+1},\dots ,\mu_n)$. Our hypothesis $A_k$ then implies that
the projectivized monodromy near a point of the stratum is finite. Equation 
\eqref{eq:normal}
shows that in the transversal direction (the $\eps$ coordinate) the multivaluedness is like that of $(\eps)^{1-|\mu'|}$. Since $\mu_i+\mu_j\in\QQ$ for all $0\le i<j\le k$ and the sum of these numbers is
$\half k(k+1)|\mu'|$, it follows that $|\mu'|\in\QQ$. So we have also finite order monodromy along the exceptional divisor. This implies that we have finite local monodromy at a point of the stratum: the stratum is elliptic. We proved in fact slightly more, namely that this local monodromy group is the one associated to the Lauricella system of type $\mu'$. So we may then invoke $A_k$ to conclude that $\widetilde{V_n^{\st}}$ is in fact smooth over this stratum.

In the general case, with a stratum corresponding to several
clusters forming, we have topologically a product situation: the local monodromy group near a point of that stratum decomposes as a product with each factor corresponding to a cluster being formed. It is clear that if each cluster is of elliptic type, then so is the stratum. Its preimage in $\widetilde{V_n^{\st}}$ will be smooth.

The asserted regularity properties of this extension of the Lauricalla map hold on codimension strata by Lemma \ref{lemma:ellhyperplane}. But then they hold everywhere, because the locus where a homolomorphic map between complex manifolds of the same dimension fails to be a local isomorphism is of codimension $\le 1$.
\end{proof}

\begin{proof}[Proof that $B_{n-1}$ implies $A_n$]
Since $B_{n-1}$ holds, it follows  that
$V^f_n$ contains $V_n-\{ 0\}$. Thus $\PP F:\PP (\widetilde{V_n})\to \PP_{n-1}$ 
is defined. The latter is a $\G$-equivariant local isomorphism with $\G$ acting on 
$\PP (\widetilde{V_n})$ with compact fundamental domain (for its orbit space is the compact $\PP(V_n)$) and on the range as a group of isometries. This implies that  $\PP (\widetilde{V_n})$ is complete. According to Lemma
\ref{lemma:covering}, $\PP F$ is then an isomorphism. 
Hence  $F:\widetilde{V_n\!\! -\!\!\{ 0\}}\to \CC^n-\{ 0\}$ is a covering projection. 
Since the domain of the latter is connected and the range is simply connected, this map is an isomorphism. In particular, $\PP (\widetilde{V_n})$ is compact, so that the covering 
$\PP (\widetilde{V_n})\to \PP(V_n)$ is finite. This  means  that
 the projectivization of $\G$ is finite. On the other hand, the $\CC^\times$-action on
$V_n-\{ 0\}$ needs a finite cover (of degree equal to the denominator of $1-|\mu|$) 
to lift to $\widetilde{V_n\!\! -\!\!\{ 0\}}$. This implies that $\G$ is finite, so that $V_n^f=V_n$.
It is now clear that $F:\widetilde{V_n}\to \CC^n$ is an isomorphism. It is $\G$-equivariant and   drops to an isomorphism $V_n\to \G\bs\CC^n$ of affine varieties.
\end{proof}

In the parabolic case $\PP(V_n)$ acquires the structure of an
orbifold modeled on flat space:

\begin{corollary}[Parabolic case]\label{cor:par}
Suppose that $|\mu |=1$ and that for all $0\le k<l\le n$,  $1-\mu_k-\mu_l$ is the reciprocal of an integer. Then $\G$ acts as a complex Bieberbach group in $\AA_{n-1}$, $V_n^f=V_n-\{ 0\}$ and $\PP F:\PP( \widetilde{V_n})\to\AA_{n-1}$ is a $\G$-equivariant isomorphism which drops to an isomorphism $\PP(V_n)\to \G\bs\AA_{n-1}$.
\end{corollary}
\begin{proof}
It follows from Proposition \ref{prop:elliptic} 
that $V^f_n$ contains $V_n-\{ 0\}$ so that $\PP F:\PP (\widetilde{V_n})\to \AA_{n-1}$ 
is defined. The latter is a $\G$-equivariant local isomorphism with $\G$ acting on the 
$\PP (\widetilde{V_n})$ with compact fundamental domain and on the range as a group of isometries. Hence $\PP (\widetilde{V_n})$ is complete. It the follows from  Lemma \ref{lemma:covering} that $\PP F$ is a $\G$-equivariant isomorphism. It also follows that
$\G$ acts on $\AA_{n-1}$ discretely with compact fundamental domain. This group
is generated by complex reflections, in particular it is a complex Bieberbach group.
Clearly,  $\PP F$ induces an isomorphism $\PP (V_n)\cong \G\bs \AA_{n-1}$.  
\end{proof}

We have also partial converses of Theorem \ref{thm:ell} and Corollary \ref{cor:par}.
They will be consequences of

\begin{lemma}\label{lemma:drop}
The Lauricella map extends holomorphically over any stable stratum contained
in $V_n^f$. 
\end{lemma}
\begin{proof}
Let $S\subset\{0,\dots ,n\}$ define an stable stratum $D_S$ (i.e., $S$ has at least two members and $\sum_{k\in S}\mu_k<1$) and assume that $D_S\subset V_n^f$.
If $0\le k<l\le n$ is contained in $S$, then $\mu_k+\mu_l\le |\mu|<1$
and so the associated monodromy transformation $T$ is according to  Corollary \ref{cor:dehnmon} a reflection over an angle $2\pi (\mu_k+\mu_l)$. Since $D_S\subset V^f_n$, we must have $\mu_k+\mu_l\in\QQ$. Lemma \ref{lemma:ellhyperplane} tells us that $F$ then extends holomorphically  over the preimage of $D_{k,l}$. The usual codimension argument implies that this is then also so over the preimage of $D_S$.
\end{proof}

\begin{proposition}\label{prop:drop}
If $|\mu|<1$ and $\G$ is finite, then the Lauricella map drops to a finite map $V_n\to \G\bs \CC^n$. 

If $|\mu|=1$, $n>1$ and $\G$ acts on the complex Euclidean space $\AA_{n-1}$ as a complex Bieberbach group, then $V_n^f=V_n$ and the Schwarz map drops to a finite map $\PP(V_n)\to \G\bs\AA_{n-1}$. 
\end{proposition}
\begin{proof}
In the elliptic case, it follows from Lemma \ref{lemma:drop} that the map $F$
drops to a map $V_n\to\G\bs\CC^n$ which exists in the complex-analytic category.   The map in question is homogeneous (relative to the natural $\CC^\times$-actions) and the preimage of $0$ is $0$. Hence it must be a finite morphism. In the parabolic case, the lemma  implies that the Schwarz map determines a map $\PP(V_n)\to \G\bs\AA_{n-1}$ which lives in the complex-analytic category. This map will be finite, because its fibers are discrete and its domain is compact.
\end{proof}

\section{The hyperbolic case}\label{sect:hyp}
Throughout this section we always suppose that $1\le |\mu|<2$. 

\subsection{A projective set-up}\label{subsect:projsetup}
An important difference with the elliptic and the parabolic cases is that $z_{n+1}=\infty$ is now of the same nature as the 
finite singular points, since we have $\mu_{n+1}=1-|\mu|\in (0,1)$. This tells us that we should treat all the points $z_0,\dots ,z_{n+1}$ on the same footing. In more precise terms, instead of taking $z_{n+1}=\infty$ and study the transformation behavior of the Lauricella integrals under the affine group $\CC^\times\ltimes \CC$ of $\CC$, 
we should let $z_0,\dots ,z_{n+1}$ be distinct, but otherwise
arbitrary points of $\PP^1$ and let the group $\PGL (2,\CC)$ take role of the affine group. This means in practice that we will sometimes allow some finite $z_k$ to coalesce with $z_{n+1}$ (that is, to fly off to infinity). For this we proceed as follows. 
Let $Z_0,\dots ,Z_{n+1}$ be nonzero 
linear forms on $\CC^2$ defining distinct points $z_0,\dots ,z_{n+1}$ of $\PP^1$.
Consider the multivalued $2$-form on $\CC^2$ defined by
\[
Z_0(\zeta)^{-\mu_0}\cdots Z_{n+1}(\zeta)^{-\mu_{n+1}}d\zeta_0\wedge d\zeta_1.
\]
Let us see how this transforms under the group $\GL(2,\CC)$. 
The subgroup $\SL(2,\CC)$ leaves $d\zeta_0\wedge d\zeta_1$ 
invariant, and so it simply transforms under $\SL(2,\CC)$ 
via the latter's diagonal action on the $(\CC^2)^{n+2}$ 
(the space that contains $Z=(Z_0,\dots ,Z_{n+1})$).
The subgroup of scalars, $\CC^\times\subset\GL(2,\CC)$ leaves the $2$-form invariant. So the form has a pole of order one at the projective line $\PP^1$ at infinity. We denote the residue of that form on $\PP^1$ by $\eta_Z$.
It is now clear, that a Lauricella function $\int_\gamma \eta_Z$ will be 
$\GL(2,\CC)$-invariant. Since the $2$-form (and hence $\eta_Z$) is 
homogeneous of degree $-\mu_k$ in $Z_k$, it follows that the 
quotient of two Lauricella functions will only depend on the  
$\GL(2,\CC)$-orbit of $(z_0,\dots ,z_{n+1})$. 

Let $\Qcal_\mu^\circ$ denote the $\SL(2,\CC)$-orbit space of the subset of 
$(\PP^1)^{n+2}$ parameterizing distinct $(n+2)$-tuples in $\PP^1$. 
This is in a natural way a smooth algebraic variety which can be 
identified with $\PP(V^\circ_n)$ (every orbit is represented by an 
$(n+2)$-tuple of which the last point is $\infty$). So we have a $\G$-covering  
$\widetilde{\Qcal}^\circ_{\mu}\to \Qcal^\circ_{\mu}$ and a local isomorphism 
$\PP F: \widetilde{\Qcal}^\circ_{\mu}\to \BB_{n-1}$.
Thus far our treatment of $z_{n+1}$  as one of the other $z_i$'s has not accomplished 
anything, but it will matter when we seek to extend it as a ramified 
covering. 

We say that $z=(z_0,\dots ,z_{n+1})\in (\PP^1)^{n+2}$ is \emph{$\mu$-stable} resp.\
 \emph{$\mu$-semistable} if the $\RR$-divisor 
 $\Div (z):=\sum_{k=0}^{n+1} \mu_k(z_k)$ has no point of weight $\ge 1$ resp.\ $>1$. Let us denote the corresponding
(Zariski open) subsets of $(\PP^1)^{n+2}$ by $U_{\mu}^\st$ resp.\ 
$U_{\mu}^\sst$. 
Notice that when $z$ is $\mu$-stable, the support of $\sum_{k=0}^{n+1} \mu_k(z_k)$ has at least three points.  This implies that the $\SL (2,\CC)$-orbit space (denoted $\Qcal_\mu^\st$) of $U_{\mu}^\st$ is in a natural manner a nonsingular algebraic variety: given  a  $\mu$-stable point $z$, we can  always pick three pairwise distinct components for use as an affine coordinate for $\PP^1$.
By means of this coordinate we get a nonempty Zariski-open subset in 
$(\PP^1)^{n-1}$ which maps bijectively to an open subset of $\Qcal_\mu^\st$. These bijections define an atlas for the claimed structure. In the semistable case, we can choose a coordinate for $\PP^1$ such that $\infty$ has weight $1$. 

Geometric Invariant Theory tells us that in case the $\mu_k$'s are all rational, one can compactify  $\Qcal_\mu^\st$ to a projective variety by adding just finitely many points: one point for each orbit containing a point whose associated divisor is $(0)+(\infty)$ or equivalently, for each splitting of $\{0,\dots ,n+1\}$ into two subsets,
each of which of total $\mu$-weight $1$. (So if no such splitting exists, then $\Qcal_\mu^\st$ is already
projective variety.) Let us denote that projective compactification by 
$\Qcal^{\sst}_{\mu}$. This is in fact a quotient of a $U_{\mu}^\sst$
with the property that each fiber is the closure of a $\SL (2,\CC)$-orbit and contains a unique closed $\SL (2,\CC)$-orbit (in the strictly 
stable case  the latter is represented by a $z$ whose divisor is   $(0)+(\infty)$).

\begin{theorem}\label{thm:hyp}
Assume that for every pair $0\le k<l\le n+1$ for which $\mu_k+\mu_l<1$,
$1-\mu_k-\mu_l$ is the reciprocal of an integer. Then the monodromy covering $\widetilde{\Qcal}^\circ_{\mu}\to \Qcal^\circ_{\mu}$ extends to 
a ramified covering $\widetilde{\Qcal}^{\st}_{\mu}\to \Qcal^{\st}_{\mu}$
and $F$ extends to a $\G$-equivariant isomorphism 
$\widetilde{\Qcal}^{\st}_{\mu}\to \BB_{n-1}$. Moreover $\G$ acts in $\BB^m$ discretely and with finite covolume; this action is with compact fundamental domain
if and only no subsequence of $\mu$ has weight $1$.
\end{theorem}

\begin{remarks}
Our hypotheses imply that the $\mu_k$'s are all rational so that
the GIT compactification $\Qcal^{\sst}_{\mu}$ makes sense.
The compactication of $\G\bs \BB_{n-1}$  that results by 
$\G\bs \BB_{n-1}\cong \Qcal^{\st}_{\mu}\subset\Qcal^{\sst}_{\mu}$ 
coincides with the \emph{Baily-Borel  compactification} of $\G\bs \BB_{n-1}$.

The  cohomology and intersection homology of the variety $\Qcal^{\sst}_{\mu}$
has been investigated by Kirwan-Lee-Weintraub \cite{klw}.
\end{remarks}

Before we begin the proof of Theorem \ref{thm:hyp} we need to know a little bit about the behavior of the complex hyperbolic metric on a complex ball near a cusp. Let $W$ be a finite dimensional complex vector space equipped with a nondegenerate Hermitian form $H$ of hyperbolic signature so that $H(w,w)>0$ defines a complex ball $\BB (W)\subset\PP(W)$. Let $e\in W$ be a nonzero isotropic vector. Since its orthogonal complement  is negative semidefinite,
every positive definite line will meet the  affine hyperplane  in $W$ defined by $H(w,e)=-1$. In this way we find an open subset $\Omega_e$ in this hyperplane which maps isomorphically onto $\BB(W)$.  This is what is called  a realization of $\BB(W)$ as a Siegel domain of the second kind. 

\begin{lemma}
The subset $\Omega$ of the affine space $H(w,e)=-1$ defined by $H(w,w)<0$ is invariant under translation by $\RR_{\ge 0} e$. If $K\subset\Omega$ is compact and measurable, then $K+\ii\RR_{\ge 0}e$ is as asubset of $\Omega$ complete and of finite volume.
\end{lemma}
\begin{proof}
This is well-known, but we outline the proof anyway. Write $e_0$ for $e$ and let  $e_1\in W$ be another isotropic vector such that $H(e_0,e_1)=1$ and denote by $W'$ the orthogonal complement of the span of  $e_0,e_1$. So if we write $w=w_0e_0+w_1e_1+w'$ with $w'\in W'$, then $\Omega$ is defined by $w_1=-1$ and $\re (w_0)>\half H(w',w')$. This shows in particular that $\Omega$ is invariant under translation by $\tau e$, when $\re (\tau)\ge 0$.
Let $K_o\subset \Omega$ be compact ball and suppose that $w\in K_o\mapsto H(w,e_1)$ is constant. If $R>0$, then the map 
$(w,y,x)\in K_o\times [-R,R]\times \RR_{\ge 0}\to (w+(x+\ii  y)e_0\in \Omega$,
is an embedding. It is straightforward to verify that the pull-back of the metric of $\Omega\cong\BB(W)$ is comparible to the `warped metric' $x^{-1}(g_\Omega |_{K_o})+ x^{-2}(dx^2+dy^2)$. From this it easily follows that 
$K_o\times [-R,R]\times \RR_{\ge 0}$ is complete and of finite volume.
Since any compact measurable $K\subset \Omega$ is covered by the image of finitely many maps $K_o\times [-R,R]\times \RR_{\ge 0}\to\Omega$ as above, the lemma
follows.
\end{proof}

It follows from Proposition \ref{prop:elliptic} that $\Qcal_{\mu}^\st\subset\Qcal_\mu^f$ and that the  Schwarz map  $\PP F:\widetilde{\Qcal}^{\st}_{\mu}\to \BB_{n-1}$ is a local isomorphism.  So $\Qcal^{\st}_{\mu}$ inherits a metric from $\BB_{n-1}$.
We need to show that $\Qcal^{\st}_{\mu}$ is complete and has finite volume.
The crucial step toward this  is:

\begin{lemma}\label{lemma:complete}
Let $0<r<n$ be such that $\mu_0+\cdots +\mu_r=1$. Denote by
$D$ the set of $(z_0,\dots ,z_n)\in\CC^{n+1}$ satisfying $|z_0|<\cdots <|z_r|<1<2<|z_{r+1}|<\cdots <|z_n|$ and $z_0+\cdots +z_r=0$.
Then $D$ embeds in $\Qcal^{\circ}_\mu$ and its closure in $\Qcal^{\st}_\mu$ is complete and of bounded volume.
\end{lemma}
\begin{proof}
That $D$ embeds in $\Qcal^{\circ}_\mu$ is clear. Let $D'\subset D$ be the open-dense subset of $z\in D$ for which 
 $|\arg (z_k)|<\pi$ for all $k$. There is a natural isotopy class of  $L$-arc systems $\delta$ for every $z\in D'$ characterized by the property that $\delta_k$ never crosses the negative real axis and $|\delta_k|$ is monotonous. This defines a lift $\tilde D'$ of $D'$ to $\widetilde{\Qcal}^{\st}_{\mu}$ so that is defined $F: \tilde D\to\CC^n$. 
For $t> 0$, denote by $D'(t)$  the set of $z\in D$ for which
$|z_0\cdots z_r|/|z_{r+1}\cdots z_n|\ge t$.
It is easy to see that $D'(t)$ has compact closure in $\Qcal^{\st}_\mu$ and so 
the closure of its preimage $\tilde D'(t)$  in $\widetilde{\Qcal}^{\st}_{\mu}$ is compact as well.  

Since $\mu_0+\cdots +\mu_r=1$, Lemma \ref{lemma:parlimit} will apply here. As in that lemma, we put   $\hat F:=z_{r+1}^{-\mu_{r+1}}\cdots z_n^{-\mu_n}F$. 
According to that lemma we have $\sum_{k=1}^r\im (w_k)\hat F_k(z)=0$. This amounts to 
saying that $H(\hat F, e_{r+1})=-\pi$, where  $e_{r+1}$ denotes the $(r+1)$th basis vector of $\CC^n$.  (For 
$H(F,G)=\sum_{1\le j<k\le n+1} \im (w_j\bar w_k) \bar G_j F_k$
and so $H(e_{r+1},G)=\sum_{1\le j\le r} \im (w_j)\bar G_j$.) We also notice that $H(e_{r+1},e_{r+1})=0$. So  $\hat F$  maps to the Siegel domain $\Omega$
defined in Lemma \ref{lemma:complete} if we take $e:=\pi^{-1}e_{r+1}$. 
Hence the lemma will follow if we show that the image of $\tilde D$
in $\Omega$  is contained in a subset of the form $K+\RR_{>0}e_{r+1}$.
Now notice that for $0<\eps<1$, $z\mapsto z_\eps$ maps $D(t)$ onto
$D'(t\eps^{n+1})$. From Lemma \ref{lemma:parlimit} we see that the coordinates 
$\hat F_k$ stay bounded on $\tilde D'$ for all $k\not=r+1$, whereas 
$\re\hat F_{r+1}|_{D'(t)}\to\infty$ as $t\to 0$. This means that $\tilde D'$
in $\Omega$  is contained in a subset of the form $K+\RR_{>0}e_{r+1}$.
\end{proof}

\begin{proof}[Proof of Theorem \ref{thm:hyp}]
The GIT compactification 
$\Qcal^\sst_\mu$ of $\Qcal^\st_\mu$ adds a point for every
permutation $\sigma$ of $\{0,\dots,n\}$ for which 
$\mu_{\sigma(0)}+\cdots +\mu_{\sigma(r)}=1$ for some $0<r<n$. 
If $\sigma$ is such a permutation, then  we have defined an open subset  
$D_\sigma\subset\Qcal^\circ_\mu$  as  in Lemma \ref{lemma:complete} 
and according to that Lemma, the closure of  $D_\sigma$ in $\Qcal^\st_\mu$ 
is complete and of finite volume. The complement in $\Qcal^\st_\mu$ of 
the union of these closures is easily seen to be compact.  
Hence $\Qcal^\st_\mu$ is complete and of finite volume. 
The theorem now follows from Lemma \ref{lemma:covering} 
(bearing in mind that $\Qcal^\sst_\mu=\Qcal^\st_\mu$ if 
and only if no subsequence of $\mu$ has weight
$1$).
\end{proof}

\subsection{Extending the range of applicability}\label{subsect:extending}

We begin with stating a  partial converse to Theorem \ref{thm:hyp}, the hyperbolic counterpart of Proposition \ref{prop:drop}:

\begin{proposition}\label{prop:hypdrop}
Suppose that $1<|\mu|<2$, $n>1$ and $\G$ acts on  $\BB_{n-1}$ as a \emph{discrete} group. Then $\G$ has finite covolume and the Schwarz map drops to a finite morphism $\Qcal^\st_{\mu}\to \G\bs\BB_{n-1}$.
\end{proposition}
\begin{proof} 
It follows from Lemma \ref{lemma:drop} that the Schwarz map is defined over $\Qcal^\st_{\mu}$ and hence drops to a map $\Qcal^\st_{\mu}\to \G\bs\BB_{n-1}$.
It follows from Lemma \ref{lemma:complete} (by argueing as in the proof of Theorem 
\ref{thm:hyp}) that $\Qcal^\st_{\mu}$ is complete as a metric orbifold and
of finite volume. This implies that $\Qcal^\st_{\mu}\to \G\bs\BB_{n-1}$ is
a finite morphism.
\end{proof}

This immediately raises the question which weight systems $\mu$  satisfy the hypotheses of Proposition \ref{prop:hypdrop}. The first  step toward the answer was taken by  Mostow himself \cite{mostow:int}, who observed that if some of the weights $\mu_k$ coincide, then the conditions of  \eqref{thm:ell}, \eqref{cor:par} and \eqref{thm:hyp}  may be relaxed, while still ensuring that $\G$ is a discrete subgroup of the relevant Lie group. The idea is this: 
if $\Scal_\mu$ denotes the group of  permutations of $\{ 0,\dots ,n+1\}$ which preserve the weights, then we should regard the Lauricella map $F$ as being multivalued on 
$\Scal_\mu\bs V_n^\circ$, rather than on $V_n^\circ$. This can make a difference, for the monodromy cover of $\Scal_\mu\bs V_n^\circ$ need not factor through $V_n^\circ$.  
We get the following variant of Lemma \ref{lemma:ellhyperplane}

\begin{lemma}\label{lemma:ellhyperplane1}
Suppose that in Lemma \ref{lemma:ellhyperplane} we have $\mu_k=\mu_l\in\QQ-\{\half\}$. Then the Lauricella map (the Schwarz map if  $|\mu|=1$) extends over the image in 
$D_{k,l}$ in   $\Scal_\mu\bs V_n^\circ$ as a local isomorphism if and only if $\half-\mu_k$
is the reciprocal of a positive integer.  
\end{lemma}

\begin{definition}\label{def:half}
We say that $\mu $ satisfies the \emph{half integrality conditions}
if whenever  for $0\le k<l\le n+1$ we have $\mu_k+\mu_l<1$, then
$(1-\mu_k-\mu_l)^{-1}$ is  an integer or in case  $\mu_k=\mu_l$, just half an integer. 
\end{definition}

This notion is a priori weaker  than Mostow's $\Sigma$INT condition, but in the end it apparently leads to the same set of weight systems.
Now  Proposition \ref{prop:elliptic} takes the following  form. 

\begin{proposition}\label{prop:elliptic1}
If $\mu$  satisfies the half integrality conditions, then $V_n^{\st}\subset V_n^f$, $\widetilde{\Scal_\mu\bs V_n^{\st}}$ is nonsingular, and the Lauricella map extends holomorphically to $\widetilde{\Scal_\mu\bs V_n^{\st}}$.
This extension has the same regularity properties as  the map it extends:
it is a local isomorphism when we are not in the parabolic case, whereas in the parabolic case, the Schwarz map defines a local isomorphism to $\AA_{n-1}$.  
\end{proposition}

This leads to  (see \cite{mostow:int} and for the present version,  \cite{chl}):

\begin{theorem}\label{thm:main}
Suppose that $\mu$  satisfies the half integrality conditions.
\begin{enumerate}
\item[\textit{ell:}] If $|\mu|<1$, then $\G$ is a finite complex reflection group in $\GL(n,\CC)$ and $F: \widetilde{\Scal_\mu\bs V}{}_n\to \CC^n$ is a $\G$-equivariant isomorphism which drops to an isomorphism $\Scal_\mu\bs V_n\to \G\bs\CC^n$. 
\item[\textit{par:}] If $|\mu |=1$, then $\G$ acts as a complex Bieberbach group in $\AA_{n-1}$, $V_n^f=V_n-\{ 0\}$ and $\PP F:\PP( \widetilde{\Scal_\mu\bs V}{}_n)\to\AA_{n-1}$ is a $\G$-equivariant isomorphism which drops to an isomorphism $\PP(\Scal_\mu\bs V_n)\to \G\bs\AA_{n-1}$.
\item[\textit{hyp:}] If $1<|\mu |<2$, then the monodromy covering $\widetilde{\Scal_\mu\bs\Qcal}{}^\circ_{\mu}\to \Scal_\mu\bs\Qcal^\circ_{\mu}$ extends to 
a ramified covering $\widetilde{\Scal_\mu\bs\Qcal}{}^{\st}_{\mu}\to \Scal_\mu\bs\Qcal^{\st}_{\mu}$
and $F$ extends to a $\G$-equivariant isomorphism 
$\widetilde{\Scal_\mu\bs\Qcal}{}^{\st}_{\mu}\to \BB_{n-1}$. Moreover $\G$ acts discretely
in $\BB^m$ and with finite covolume. 
\end{enumerate}
\end{theorem}

\begin{example}\label{example:half}
Let us take $n\le 10$ and $\mu_k=\frac{1}{6}$ for $k=0,\dots ,n$. So we have
$\mu_{n+1}= \frac{11-n}{6}$. The half integrality conditions  are fulfilled for all
$n\le 10$ with $1\le n\le 4$, $n=5$, $6\le n\le 11$ yielding an elliptic, parabolic 
and hyperbolic case, respectively and $\Scal_\mu$ is the
permutation group of $\{0,\dots, n\}$ for $n\le 9$ and the one of 
$\{0,\dots, 11\}$ for $n=10$.  
\end{example}

Mostow  subsequently showed  that in the hyperbolic range with $n\ge3$ we thus find all but ten of the discrete monodromy groups of finite covolume: one is missed for $n=4$ (namely $(\frac{1}{12}, \frac{3}{12},
\frac{5}{12},\frac{5}{12},\frac{5}{12},\frac{5}{12})$) and nine for $n=3$ (see \cite{mostow:disc}, (5.1)). He conjectured that in these nine cases $\G$ is always commensurable with a group obtained from Theorem \ref{thm:main}. This was proved
by his student Sauter  \cite{sauter}. It is perhaps no surprise that things are a bit different  when $n=2$ (so that we are dealing with discrete groups of automorphism of the unit disk): indeed, the exceptions then make up a number of infinite series (\cite{mostow:disc}, Theorem 3.8). It turns out that for $n>10$ the monodromy 
group is never discrete  and that for $n=10$ this happens only when  $\mu_k=\frac{1}{6}$ for $k=0,\dots ,10$. (It is not known whether there exist discrete subgroups of isometry groups of  finite covolume of a complex ball of dimension $\ge 10$.)

\section{Modular interpretation}\label{sect:modular}
We assume here that we are in the $\QQ$-hyperbolic case: 
$\mu_k\in (0,1)$ and rational  for $k=0,\dots ,n+1$ (with $\sum_{k=0}^{n+1}\mu_k=2$ as always). 

\subsection{Cyclic covers of $\PP^1$}\label{subsect:cycliccover}
We will show that the Schwarz map can be interpreted as  a `fractional period' map. This comes about  by passing to a cyclic cover of $\PP^1$
on which the Lauricella integrand becomes a regular differential. Concretely,
write $\mu_k=d_k/m$ with $d_k,m$ positive integers such that the $d_k$'s
have no common divisor, and  write $m_k$ for the denominator of $\mu_k$.
Consider the cyclic cover $C\to \PP^1$ of order 
$m$ which has ramification over $z_k$ of order $m_k$. 
In affine coordinates, $C$ is given as the normalization of the curve defined 
by
\[
w^m=\prod_{k=0}^n (z_k-\zeta)^{d_k}.
\]
This is a cyclic covering which has the group $G_m$ of $m$th roots of unity
as its Galois group: $g^*(w,z)=(\chi (g)w,z)$, where 
$\chi :G_m\subset\CC^\times$ stands for the tautological character. The Lauricella integrand pulls  back to a univalued differential $\tilde\eta$ on $C$, 
represented by $w^{-1}d\zeta$ so that $g^*(\tilde\eta)=\bar\chi (g)\tilde\eta$.
Hence, if we let $G_m$ act on forms in the usual manner ($g\in G_m$ acts 
as  $(g^{-1})^*$), then $\tilde\eta$ is an eigenvector with character $\chi$.
It is easily checked that $\tilde\eta$ is regular  everywhere. 

In order to put this in a period setting, we recall some generalities concerning the Hodge decomposition of $C$: its space of holomorphic 
differentials, $\Omega (C)$, has dimension equal to the genus $g$ of $C$ and
$H^1(C;\CC)$ is canonically represented 
on the form level by the direct sum $\Omega (C)\oplus\overline{\Omega}(C)$
(complex conjugation on forms corresponds to complex conjugation in 
$H^1(C;\CC)$ with respect to $H^1(C;\RR)$). The intersection product on
$H^1(C;\ZZ)$ defined by $(\alpha,\beta)\mapsto (\alpha\cup \beta)[C]$ (where
the fundamental class $[C]\in H_2(C,\ZZ)$ is specified by the complex orientation of $C$), is on the
form level given by $\int_C \alpha\wedge\beta$. The associated Hermitian form
on $H^1(C;\CC)$ defined by $h(\alpha,\beta):=\frac{\ii}{2}
(\alpha\cup \bar\beta)[C]=\frac{\ii}{2}\int_C\alpha\wedge\beta$ has signature
$(g,g)$. The Hodge decomposition $H^1(C;\RR)=\Omega (C)\oplus
\overline{\Omega}(C)$ is $h$-orthogonal with the first summand positive 
definite and the second negative definite. The Hodge decomposition, the 
intersection product and (hence) the Hermitian form $h$ are all 
left invariant by the action of $G_m$. 
 
\begin{proposition}\label{prop:signature}
The eigenspace $\Omega(C)^\chi$ is of dimension one and spanned by 
$\tilde\eta$ and the eigenspace $\overline{\Omega}(C)^\chi$ is of dimension $n-1$. The eigenspace $H^1(C,\CC)^\chi$ has signature $(1,n-1)$ and we have 
$h(\tilde\eta,\tilde\eta)=-mN(F(z),F(z))$.
\end{proposition}

\begin{lemma}\label{lemma:eigenforms}
Let $r\in\{ 0,1,\dots ,m-1\}$. Then the  eigenspace $\Omega (C)^{\chi^r}$ is
spanned by the forms $w^{-r}f(\zeta)d\zeta$ where $f$ runs over the
polynomials of degree $<-1+r\sum_{k=0}^n\mu_k$ that have in 
$z_k$ a zero of order $\ge [r\mu_k]$, $k=0,\dots ,n$. In particular,
$\dim \Omega (C)^{\chi^r}$ is the largest integer smaller than $\sum_{k=0}^n \{r\mu_k\}$ (recall that $\{ a\}:=a-[a]$). 
\end{lemma}
\begin{proof}
Any meromorphic differential on $C$ which transforms according to the
character $\chi^r$, $r=0,1,\dots ,m-1$, is of the
form $w^{-r}f(\zeta)d\zeta$ with $f$ meromorphic. 
A local computation shows that in order that such a differential be regular, 
it is necessary and sufficient 
that $f$ be a polynomial of degree $<-1+r\sum_{k=0}^n\mu_k$ which has in 
$z_k$ a zero of order $>-1+r\mu_k$, that is, of order $\ge [r\mu_k]$. 
\end{proof}

\begin{proof}[Proof of Proposition \ref{prop:signature}]
If we apply Lemma \ref{lemma:eigenforms} to the case $r=1$, then we find that $f$ must have degree $<-1+\sum_{k=0}^n\mu_k=
1-\mu_{n+1}$ and as $\mu_{n+1}\in (0,1)$, this means that $f$ is 
constant. So $\tilde\eta$ spans $\Omega(C)^\chi$.

For $r=m-1$, we find  that $\dim\Omega(C)^{\bar\chi}$
 is the largest integer smaller than $\sum_{k=0}^n \{(m-1)\mu_k\}=
\sum_{k=0}^n \{d_k-\mu_k\}=\sum_{k=0}^n (1-\mu_k)=n-1+\mu_{n+1}$,
that is, $n-1$. Since $\overline{\Omega}(C)^\chi$ is the complex conjugate of 
$\Omega(C)^{\bar\chi}$, it follows that this space has dimension $n-1$ also.

That $H^1(C,\CC)^\chi$ has signature $(1,n-1)$ is now a consequence 
of its orthogonal decomposition into $\Omega(C)^\chi$ and 
$\overline{\Omega}(C)^\chi$. Finally,
\[
h(\tilde\eta,\tilde\eta)=\frac{\ii}{2}\int _C \tilde\eta\wedge \overline{\tilde\eta}=
\frac{m\ii}{2}\int_{\CC} \eta\wedge \bar\eta=-mN(z,z) (>0).\qedhere
\]
\end{proof}

So the Schwarz map $\PP F :\widetilde{\Qcal}_{\mu}^\st\to \BB_{n-1}$ can now be understood as attaching to the curve $C$ with its $G_m$-action the Hodge decomposition of $H^1(C;\CC)^\chi$. 

\subsection{Arithmeticity}\label{subsect:arithmetic}
The above computation leads to the following arithmeticity criterion for $\G$:

\begin{theorem}\label{thm:arithmetic}
The monodromy group $\G$ is arithmetic if and only if for every 
 $r\in(\ZZ/m)^\times-\{\pm 1\}$  we have $\sum_{k=0}^n \{ r \mu_k\} \le 1$ or $\sum_{k=0}^n \{ -r\mu_k\} \le 1$.
\end{theorem}

We need the following density lemma.

\begin{lemma}\label{lemma:dense}
The Zariski closure of $\G$ in  
$\GL(H^1(C,\CC)^{\chi}\oplus H^1(C,\CC)^{\bar\chi})$ is 
defined over $\RR$ and the image of its group of real points in  
the general linear group of $H^1(C,\CC)^{\chi}$ contains the 
special unitary group of $H^1(C,\CC)^{\chi}$.
\end{lemma}

The proof amounts to exhibiting sufficiently many complex reflections in $\G$.
It is somewhat technical and we therefore omit it.

\begin{proof}[Proof of Theorem \ref{thm:arithmetic}]
Let us abbreviate $H^1(C,\CC)^{\chi^r}$ by $H_r$.
The smallest subspace of $H^1(C,\CC)$ which contains 
$H_1$ and is defined over $\QQ$ is the sum of the eigenspaces
$H:=\oplus_{r\in(\ZZ/m)^\times}H_r$. We may identify $H$ with the quotient 
of $H^1(C,\CC)$ by the span of the images of the maps  
$H^1(G_k\bs C,\CC)\to H^1(C,\CC)$, where
$k$ runs over the divisors $\not=1$ of $m$. In particular,
$H(\ZZ):=H^1(C,\ZZ)\cap H$ spans $H$. The monodromy group $\G$ may be regarded as a subgroup of $\GL (H_\ZZ)$. 
On the other hand, $\G$ preserves each summand $H_r$. So if we denote by $\Gcal$ the $\QQ$-Zariski closure of $\G$ in $\GL (H)$, then
$\G\subset \Gcal (\ZZ)$ and $\Gcal (\CC)$ decomposes as $\Gcal(\CC)=\prod_{r\in(\ZZ/m)^\times} \Gcal_r(\CC)$ with 
$\Gcal_r(\CC)\subset\GL(H_r)$. To say that $\G$ is arithmetic is to say that $\G$ is of finite index in $\Gcal(\ZZ)$.

Since $H_r\oplus H_{-r}$ is defined over $\RR$, so is $\Gcal_{r,-r}:=\Gcal_r\times \Gcal_{-r}$. According to Lemma \ref{lemma:dense}, the image of $\Gcal_{1,-1}(\RR)$ in $\Gcal_r(\CC)$ contains the special unitary group of $H_1$. The summand  $H_r$ with its Hermitian form  is for $r\in(\ZZ/m)^\times$  a Galois conjugate of $H_1$ and so it then follows that the image of $\Gcal_{r,-r}(\RR)$ in $\Gcal_r(\CC)$ contains the  special unitary group of  $H_r$. 

Suppose now that $\G$ is arithmetic.
The projection $\Gcal(\RR)\to\Gcal_{1,-1}(\RR)$ is injective on $\G$ and so the kernel of this projection must be anisotropic: $\Gcal_{r,-r}(\RR)$ is compact  for $r\not=\pm 1$. This means that the Hermitian form on $H_r$ is definite for $r\not= \pm 1$. Since $H_r=\Omega (C)^{\chi^r}\oplus \overline{\Omega (C)^{\chi^{-r}}}$
with the first summand positive and the second summand negative, this means that
for every $r\in (\ZZ/m)^\times-\{\pm 1\}$ (at least) one of the two summands must be  trivial. Following  Lemma \ref{lemma:eigenforms} this amounts to $\sum_{k=0}^n \{ r \mu_k\} <1$ or $\sum_{k=0}^n \{-r\mu_k\} <1$.

Suppose conversely,  that for all  all $r\in(\ZZ/m)^\times-\{\pm 1\}$ we have $\sum_{k=0}^n \{ r \mu_k\} <1$ or $\sum_{k=0}^n \{ -r\mu_k\} <1$. As we have just seen, this amounts to $\Gcal_{r,-r}(\RR)$ being compact for all $r\in(\ZZ/m)^\times-\{\pm 1\}$.
In other words, the projection, $\Gcal (\RR)\to \Gcal_{1,-1}(\RR)$ has compact kernel. Since $\Gcal (\ZZ)$ is discrete in $\Gcal (\RR)$, it follows that its image in $\Gcal_{1,-1}(\RR)$ is discrete as well. In particular, $\G$ is discrete in $\GL(H_1)$. Following Proposition \ref{prop:hypdrop} this implies that $\G$ has finite covolume in $\Gcal_{1,-1}(\RR)$. Hence it has also in  finite covolume in $\Gcal(\RR)$. This implies that $\G$ has finite index in
$\Gcal (\ZZ)$.
\end{proof}

\begin{example}
The case for which $n=3$, $(\mu_0,\mu_1,\mu_2,\mu_3)=(\frac{3}{12},\frac{3}{12},\frac{3}{12},
\frac{7}{12})$ (so that $\mu_4=\frac{8}{12}$) satisfies the hypotheses of Theorem \ref{thm:hyp}, hence yields a monodromy group which operates on $\BB_2$ discretely with compact fundamental domain. But  the group is not arithmetic since 
we have both $\sum_{k=0}^3 \{ 5\mu_k\}=\frac{5}{3}>1$ and 
$\sum_{k=0}^3 \{ -5\mu_k\}=\frac{7}{3}>1$.
\end{example}

\subsection{Working over a ring of cyclotomic integers}\label{subsect:cyclotomic}

If we are given an $L$-arc system $\delta$, then $C\to \PP^1$ comes with 
a section (continuous outside $\delta$) in much the same way as  we found
a determination of $\eta_z$: for $\zeta$ in a left half plane, $\prod_{k=0}^n (z_k-\zeta)^{d_k}$ has argument $<\pi/2$ in absolute value and so has there  a natural $m$th root (with argument $<\pi/2m$ in absolute value); the resulting section we find there is then extended in the obvious way.
We identify  $\delta_k$ with its image in $C$ under the section and thus regard 
it as a chain on $C$. For $k=1,\dots ,n$, we introduce a $\ZZ[\zeta_m]$-valued $1$-chain on $C$:
\[
\eps_k:= \bar w_k\sum_{g\in G_m} \chi (g) g_*\delta_k.
\]
Notice that the coefficient $\bar w_k$ is an $m$th root of unity and so a unit of
$\ZZ[\zeta_m]$. We put it in, in order to maintain the connection with the 
Lauricella map. It will also have the effect of keeping some of the formulae simple. 

\begin{lemma}
The element $\eps_k$
is a $1$-cycle  on $C$ with values in $\ZZ[\zeta_m]$ and has the property that  $g_*\eps_k=\bar\chi (g)\eps_k$ 
(and hence defines an element of $H_1(C,\ZZ[\zeta_m])^{\bar\chi}$). 
We have $\int_{\eps_k}\tilde \eta=mF_k(z,\delta)$.  Moreover,
$H_1(C,\ZZ[\zeta_m])^{\bar\chi}$ is as a $\ZZ[\zeta_m]$-module freely generated
by $\eps_1,\dots ,\eps_n$.
\end{lemma}
\begin{proof}
The identity involving integrals  is verified by
\begin{multline*}
\int_{\eps_k}\tilde \eta=\bar w_k\sum_{g\in G_m}\chi (g) 
\int_{g_*\delta_k}\tilde \eta=
\bar w_k\sum_{g\in G_m}\chi (g) 
\int_{\delta_k}g^*\tilde \eta=\\
=\bar w_k\sum_{g\in G_m}\chi (g) \int_{\delta_k}\bar\chi (g)\eta =
m\bar w_k\int_{\delta_k}\eta= mF_k(z,\delta).
\end{multline*}
Give $\PP^1$ the structure of a finite cell complex by taking 
the singletons $\{z_0,\dots ,z_n\}$ as $0$-cells, the intervals $\delta_1,\dots ,\delta_n$ minus their end points as $1$-cells and $\PP^1-\cup_{i=k}^n\delta_k$ as $2$-cell. The connected components of the preimages of cells in $C$ give the latter the structure of a finite cell complex as well (over the $2$-cell we have one point of ramification, namely $\infty$,
and so connected components of its  preimage  are indeed $2$-cells). The
resulting cellular chain complex of $C$, 
\[
0\to C_2\to C_1\to C_0\to 0,
\]
comes with a $G_m$-action. Notice that $C_1$ is the free $\ZZ[G_m]$-module generated by $\delta_1,\dots ,\delta_n$. On the other hand,
$C_0\cong\oplus_{k=0}^n \ZZ[G_m/G_{m_k}]$ and 
$C_2\cong\ZZ[G_m/G_{m_{n+1}}]$, so that $(C_0)^{\bar\chi}=
(C_2)^{\bar\chi}=0$. The remaining assertions of the lemma follows from this.
\end{proof}

We explicitly describe the Hermitian form on the free 
$\ZZ[\zeta_m]$-module $H_1(C,\ZZ[\zeta_m])^{\bar\chi}$:

\begin{proposition}\label{prop:H}
The  Hermitian form $H=-\frac{1}{m}h$ is on the basis $(\eps_1,\dots ,\eps_n)$ given as follows: for $1\le l\le k\le n$ we have 
\[
H(\eps_k,\eps_l)=
\begin{cases}
0 &\text{ if $l<k-1$,}\\
-\frac{1}{4}\sin (\pi/m)^{-1} & \text{ if $l=k-1$,}\\
\frac{1}{4}(\cot (\pi/m_{k-1})+\cot (\pi/m_{k})) & \text{ if $l=k$.}
\end{cases}
\]
\end{proposition} 

It is perhaps noteworthy that this proposition shows that the matrix  
of $H$ on $\eps_1,\dots ,\eps_n$ only involves the denominators of the weigths $\mu_0,\dots ,\mu_n$.
The proof relies on a local computation of intersection multiplicities with
values in $\ZZ[\zeta_m]$. The basic situation is the following. Consider the 
$G_m$-covering $X$ over the complex unit disk $\Delta$ defined by $w^m=z^d$, where $d\in\{ 1,\dots ,m-1\}$ and $g\in G_m$ acts as $g^*w=\chi (g)w$.
The normalization $\tilde X$ of $X$ consists of $e:=\gcd (d,m)$ copies $\Delta$,
$\{\Delta_k\}_{k\in \ZZ/e}$,  as follows: if we write $m=e\bar m$ and $d=e\bar d$ and $t_k$ is the coordinate of $\Delta_k$, then $\Delta_k\to X$ is given by 
$z=t_k^{\bar m}$ and $w=\zeta_m^kt_k^{\bar d}$, so that on $\Delta_k$, 
$w^{\bar m}=\zeta_m^{k\bar m}t_k^{\bar d\bar m}=\zeta_e^{k}z^{\bar d}$. 
If $g_1\in G_m$ is such that $\chi (g_1)=\zeta_m$, then $g_1^*(t_{k+1})=t_k$ 
$k=0,1\dots ,e-1$ and $g_1^*t_0=\zeta_mt_{e-1}$ (because 
$w|\Delta_{k+1}=\zeta_m^{k+1}t_{k+1}^{\bar d}$ and $(g_1^*w)|\Delta_{k}=\zeta_m w|\Delta_k=\zeta_m^{k+1}t_k^{\bar d}$).

Choose $\theta\in (0,2\pi)$ and let $\delta$ resp.\ $\delta'$ be the ray on $\Delta_0$ defined by 
$t_0=r$ resp.\ $t_0=r\exp (\ii\theta /\bar m)$ with $0\le r <1$.
We regard either as a chain with closed support. 
Notice that $z$ maps $\delta$ resp.\ $\delta'$ onto $[0,1)$ resp. a ray $\not= [0,1)$.
Consider the $\ZZ[\zeta_m]$-valued chains with closed support
\[
\tilde\delta :=\sum_{g\in G_m} \chi (g)g_*\delta, \quad
\tilde\delta' :=\sum_{g\in G_m} \chi (g)g_*\delta'.
\]
These are in fact $1$-cycles with closed support which only meet
in the preimage of the origin (a finite set). So they have a well-defined intersection number.

\begin{lemma}\label{lemma:H}
We have $\tilde\delta\cdot \overline{\tilde\delta'}= m\zeta_m(\zeta_m-1)^{-1}=
\half m (1-\ii\cot(\pi/\bar m))$.
\end{lemma}
\begin{proof}
This intersection product gets a contribution from each connected component $\Delta_k$. Because of the $G_m$-equivariance these contributions are the
same and so it is enough to show that the contribution coming from one of 
them is $(m/2e) (1+\ii\cot(\pi/2\bar m))=\half \bar{m} (1+\ii\cot(\pi/2\bar m))$. This means that there is no loss in generality in assuming that $d$ and $m$ are relative prime. Assuming that this is the case, then we can compute the intersection product if  
we write $\tilde\delta$  and  $\tilde\delta'$ as a sum of
closed $1$-cycles with coefficients in $\ZZ[\zeta_m]$. This is accomplished by
\begin{multline*}
\tilde\delta =\sum_{g\in G_m} \chi (g)g_*\delta=\\
=\sum_{k=1}^m (1+\zeta_m+\cdots +\zeta_m^{k-1}) (g_{1*}^{k-1}\delta-g_{1*}^{k}\delta)=
\sum_{k=1}^{m} \frac{1-\zeta_m^k}{1-\zeta_m} (g_{1*}^{k-1}\delta-g_{1*}^{k}\delta),
\end{multline*}
(notice that $g_{1*}^{k-1}\delta-g_{1*}^{k}\delta$ is closed, indeed)
and likewise for $\tilde\delta'$. We thus reduce our task to computing the intersection numbers $(g_{1*}^{k-1}\delta-g_{1*}^{k}\delta)\cdot(g_{1*}^{l-1}\delta'-g_{1*}^{l}\delta')$. This is easy: we find that this equals $1$ if $l=k$, $-1$ if $l=k-1$
and $0$ otherwise. Thus
\begin{multline*}
\tilde\delta\cdot \overline{\tilde\delta'}=
\sum_{k=1}^{m} \frac{1-\zeta_m^k}{1-\zeta_m} \bar\zeta_m^{k-1}=
\frac{m\zeta_m}{\zeta_m-1}=\frac{m\zeta_{2m}}{\zeta_{2m}-\bar\zeta_{2m}}=\half m (1-\ii\cot(\pi/m)). \qedhere
\end{multline*}
\end{proof}

\begin{proof}[Proof of \ref{prop:H}]
We may of course assume that each $z_k$ is real: $z_k=x_k\in\RR$ with
with $x_0<x_1<\cdots <x_n$ and that $\delta_k=[x_{k-1},x_k]$. 
Let  us  put $\tilde\delta_k:=w_k\eps_k=\sum_{g\in G_m}\chi (g)g_*\delta_k$ and 
compute $\tilde\delta_k\cdot\overline{\tilde\delta_l}$ for 
$1\le l\le k\le n$.  It is clear that this is zero in case $l<k-1$.
For $l=k$, we let $\delta'_k$ go in a straight line from $x_{k-1}$ to a point in the upper half plane (with real part $\half x_{k-1}+\half x_k$, say) and then straight to $x_k$. We have a naturally defined $\ZZ[\zeta_m]$-valued
$1$-chain $\tilde\delta'_k$ on $C$ homologous to $\tilde\delta_k$ and with support lying over $\delta_k$. So $\tilde\delta_k\cdot\overline{\tilde\delta_k}=
\tilde\delta_k\cdot\overline{\tilde\delta'_k}$. The latter is computed with the help
of Lemma \ref{lemma:H}: the contribution over $x_{k-1}$ is 
$\half m (1-\ii\cot(\pi/m_{k-1}))$ and over $x_k$ it is  $-\half m (1-\ii\cot(\pi/m_k))$
and so $\eps_k\cdot\eps_k=\tilde\delta_k\cdot\overline{\tilde\delta'_k}=
-\half m \ii\cot(\pi/m_{k-1}))+\half m \ii\cot(\pi/m_k)$. 
We now do the case $l=k-1$. The $1$-chains on $C$ given by  
$\delta_{k-1}$ and $\delta_k$ make an angle over $x_{k-1}$ of 
$\pi\mu_{k-1}=\pi d_{k-1}/m$. In terms of the local picture of 
Lemma \ref{lemma:H} this means that the pair $(\delta_k,\delta_{k-1})$ 
corresponds to $(\delta ,-\bar\zeta_{2m}^{d_{k-1}-1}\delta')$. It follows that 
\begin{multline*}
\tilde\delta_k\cdot\overline{\tilde\delta}_{k-1}=
\tilde\delta\cdot\overline{-\bar\zeta_{2m}^{d_{k-1}-1}\tilde\delta'}=
-\zeta_{2m}^{d_{k-1}-1}\tilde\delta\cdot\overline{\tilde\delta'}=\\
=-\zeta_{2m}^{d_{k-1}-1}m\zeta_m(\zeta_m-1)^{-1}
=-m(\zeta_{2m}-\bar\zeta_{2m})^{-1}e^{\ii\pi\mu_{k-1}}.
\end{multline*}
Hence $\eps_k\cdot\overline\eps_{k-1}=-m(\zeta_{2m}-\bar\zeta_{2m})^{-1}$ and so 
$H(\eps_k,\eps_{k-1})=-\frac{1}{2m\ii}\eps_k\cdot\overline\eps_{k-1}=
(2\ii (\zeta_{2m}-\bar\zeta_{2m}))^{-1}=-\frac{1}{4}(\sin(\pi/m))^{-1}$ is as asserted.
\end{proof}

\section{Generalizations and other view points}\label{sect:other}

\subsection{Higher dimensional integrals}\label{subsect:higherdim}

This refers to the situation where $\PP^1$ and the subset
$\{ z_0,\dots ,z_{n+1}\}$ are replaced by a projective arrangement;
such generalizations were considered by Deligne, Varchenko \cite{var} and others. 
To be specific, fix an integer $N\ge 1$, a finite set $K$ with at least $N+2$ elements
and a \emph{weight function} $\mu :k\in K\mapsto \mu_k\in(0,1)$.
Given an injective map $z:k\in K\mapsto z_k\in\check{\PP}^N$, choose 
for every $k\in K$ a linear form  $Z_k:\CC^{N+1}\to \CC$ whose zero set is  the hyperplane $H_{z_k}$ defined by $z_k$ and put 
\[
\eta_z =\res_{\PP^N} \left(\prod_{k\in K} Z_k(\zeta)^{-\mu_k}\right)d\zeta_0\wedge\cdots\wedge d\zeta_N.
\]
This is a multivalued holomorphic $N$-form on $U_z:=\PP^N-\cup_{k\in K} H_{z_k}$. If $\sigma$ is
a sufficiently regular relative $N$-chain of the pair $(\PP^N, \PP^N-U_z)$
and we are given a determination of $\eta$ over $\sigma$, then $\eta$ is integrable 
over $\sigma$ so that $\int_\sigma\eta$ is defined. Here it pays however to take the more cohomological approach that we briefly described in Subsection \ref{subsect:localsystem}. So we 
let $\LL_z$ be the rank one local system on $U_z$ such that 
its monodromy around $H_{z_k}$ is multiplication by $\exp (2\pi\mu_k\ii)$ and 
endow it with a flat Hermitian metric.
Then after the choice of a multivalued section of $\LL_z$ of unit norm, 
$\eta_z$ can be interpreted as a section of $\Omega^N_{U_z}\otimes_\CC\LL_z$. 
It thus determines an element $[\eta_z]\in H^N(\LL_z)$. Similarly,  
$\sigma$ plus the determination of $\eta_z$ over $\sigma$ defines an element  
$[\sigma]\in H_N(\PP^N,\PP^N-U_z;\LL_z^\vee)$. The latter space is dual to 
$H^N(\LL_z)$ by Alexander duality in such a manner that $\int_\sigma\eta_z$ 
is the value of the Alexander pairing on $([\eta_z],[\sigma])$. 
In order that $\eta_z$ is square integrable it is necessary and s
ufficient that for every nonempty  intersection $L$ of hyperplanes 
$H_{z_k}$ we have
$\sum_{\{ k\, |\, H_{z_k}\supset L\}} \mu_k<\codim (L)$. Assume  
that this is the case. Then $\eta_z$ defines in fact a class in the 
intersection homology space 
$\IH ^m(\PP^N,\LL_z)$. This space comes a natural hermitian form $h$ 
for which  $h(\eta_z,\eta_z)>0$. (It is clear that the line spanned 
by $\eta_z$ only depends $z$; Hodge theory tells us that  
the image of that line is $F^NIH^N(\PP^N,\LL)$.) 
So in order that the situation is like the one we studied we would 
want that the orthogonal complement of $\eta_z$ in  
$\IH^N(\PP^N,\LL_z)$ to be negative. Unfortunately this seems rarely to 
be the case when $N>1$. When that is so, then we might vary 
$z$ over the connected constructible set $S$ of injective maps 
$K\to \check{\PP}^N$ for which the the topological type of the  
arrangement it defines stays constant.  Then over $S$  we have a local 
system $\HH_S$ whose stalk at $z\in S$  is $\IH^N(\PP^N,\LL_z)$ and the 
Schwarz map which assigns to $z$ the line in $\HH_z$ defined by $\eta_z$ 
will take values in a ball. The first order of business should be to 
determine the cases for which the associated monodromy group is discrete. 

\subsection{Geometric structures on arrangement complements}\label{subsect:geomstr}
In \cite{chl} Couwenberg, Heckman and I developed a generalization  of the Deligne-Mostow theory that starts with a slightly different point of view. The point of departure is here a finite dimensional complex inner product 
space $V$, a finite collection $\Hcal$ of linear hyperplanes in $V$ and a map
$\kappa$ which assigns to every $H\in\Hcal$ a positive real number $\kappa_H$. 
These data define a connection $\nabla^\kappa$ on the tangent bundle of the 
arrangement complement $V^\circ:=V-\cup_{h\in\Hcal}H$ as follows. For $H\in\Hcal$ 
denote by $\pi_H\in\End (V)$ the orthogonal projection with kernel $H$ 
and by $\omega_H$ the logarithmic differential on $V$
defined  by $\phi_H^{-1}d\phi_H$, where $\phi_H$ is a linear form on $V$ with 
kernel $H$. Form $\Omega^\kappa:= 
\sum_{H\in\Hcal} \kappa_H\pi_H\otimes \omega_H$ and regard it as a differential on 
$V^\circ$ which takes values in the tangent bundle of $V^\circ$, or rather, as  a connection form on this tangent bundle: a connection is defined by
\[
\nabla^\kappa:=\nabla^0-\Omega^\kappa,
\]
where $\nabla^0$ stands for the usual affine connection on $V$ restricted to 
$V^\circ$. This connection is easily verified to be torsion free. It is 
well-known that such a connection defines an affine structure (that is, it 
defines an atlas of charts whose transition maps are affine-linear) precisely 
whenthe connection is flat; the sheaf of affine-linear functions are then the 
holomorphic functions whose differential is flat for the connection 
(conversely, an affine structure is always given by a flat torsion free 
connection on the tangent bundle). There is a simple criterion  for the 
flatness of $\nabla^\kappa$ in terms 
of linear algebra. Let $\Lcal(\Hcal)$ denote the collection of 
subspaces of $V$ that are intersections of members of $\Hcal$
and let for  $L\in\Lcal(\Hcal)$ $\Hcal_L$ be the set of $H\in\Hcal$ 
containing $L$. Then the following properties are equivalent: 
\begin{enumerate}
\item[(i)] $\nabla$ is flat,
\item[(ii)] $\Omega\wedge \Omega=0$,
\item[(iii)]  for every pair  $L,M\in \Lcal(\Hcal)$ with $L\subset M$, 
the endomorphisms $\sum_{H\in\Hcal_L} \kappa_H\pi_H$ and
$\sum_{H\in\Hcal_M}  \kappa_H\pi_H$ commute,
\item[(iv)] for every $L\in \Lcal(\Hcal)$ of codimension $2$, the sum
$\sum_{H\in\Hcal_L} \kappa_H\pi_H$ commutes with each of its terms.
\end{enumerate}
If these mutually equivalent conditions are satisfied we call the triple $(V,\Hcal,\kappa)$ a \emph{Dunkl system}. 

Suppose that $(V,\Hcal,\kappa)$ is such a system so that $V^\circ$ comes with an affine
structure. If $L\in \Lcal(\Hcal)$ is irreducible
(in the sense that there is no nontrivial decomposition of $\Hcal_L$ such that
the corresponding intersections are perpendicular), then the fact that
$\sum_{H\in\Hcal_L} \kappa_H\pi_H$ commutes with each of its terms implies
that this sum must be proportional to the orthogonal projection with kernel
$L$, $\pi_L$. A trace computation shows that the sclalar factor must be
$\kappa_L:=\codim (L)^{-1}\sum_{H\in\Hcal_L}\kappa_H$. 
Let us now assume that the whole system is irreducible in the sense that
the intersection of all members of $\Hcal$ is reduced to the origin and
that this intersection is irreducible. We then have defined $\kappa_0=
\dim (V)^{-1}\sum_{H\in\Hcal}\kappa_H$. 
The connection is invariant under scalar multiplication by $e^t\in\CC^\times$
and one verifies that for $t$ close to $0$, the corresponding affine-linear
transformation is like scalar multiplication by $e^{(1-\kappa_0)t}$ if
$\kappa_0\not=1$ and by a translation if $\kappa_0=1$. This means that 
if $\kappa_0\not=1$, the affine structure on $V^\circ$ is in fact a linear structure
and that this determines  a (new) projective structure on $\PP(V^\circ)$, whereas
when $\kappa_0=1$ (the \emph{parabolic} case), $\PP(V^\circ)$ inherits
an affine structure which makes the projection $V^\circ\to \PP(V^\circ)$ affine-linear.
Notice that if $(V,\Hcal,t\kappa)$ will be a Dunkl system for every $t>0$.  The behavior of that
system (such as its monodromy) may change dramatically if we vary $t$.
 
Before we proceed, let us show how a weight system $\mu$
that gives rise to the Lauricella differential also gives rise to such an
irreducible Dunkl system:
we take $V=V_n=\CC^{n+1}/\text{main diagonal}$, $\Hcal$ will be the
collection of diagonal hyperplanes $H_{k,l}:=(z_k=z_l)$, $0\le k<l\le n$, and $\kappa (H_{k,l})=\mu_k+\mu_l$. The inner product on $V_n$ comes from the inner product
on $\CC^{n+1}$ for which $\la e_k,e_l\ra =\mu_k\delta_{k,l}$ and is the one which
makes the projection $\CC^{n+1}\to V_n$ selfadjoint. It is an amusing exercise to
verify that the connection is flat indeed and that the space of affine-linear 
functions at $z\in V_n^\circ$ is precisely the space of solutions of the 
system of differential equations we encountered in part (c) of 
Proposition \ref{prop:elementary}. So the Schwarz map
is now understood as a multivalued chart  (in standard terminology, a developing map) for the new projective structure on $\PP(V^\circ_n)$. We also find that $\kappa_0=|\mu|$; more generally, an irreducible member $L\in L(\Hcal)$ is given by a subset
$I\subset\{0,\dots, n\}$ with at least two elements 
(so that $L=L(I)$ is the locus where all $z_k$, $k\in I$ coincide) and 
$\kappa_{L(I)}=\sum_{k\in I}\mu_k$. 

Another interesting class of examples is provided by the finite complex reflection 
groups: let $G$ be a finite complex reflection group operating irreducibly and unitarily in 
a complex inner product space $V$, $\Hcal$ the collection of 
complex hyperplanes of $G$ and $H\in\Hcal\mapsto\kappa_H$ constant on the
$G$-orbits. Then $(V,\Hcal,\kappa)$ is a Dunkl system. 

It turns out that in many cases of interest (including the examples mentioned above),
one can show that there exists   a $\nabla^\kappa$-flat Hermitian form $h$ on $V^\circ$
with the following properties
\begin{description}
\item[ell] if $0<\kappa_0<1$, then $h$ is positive definite, 
\item[par] if $\kappa_0=1$, then $h$ positive semidefinite with kernel the tangent spaces to the $\CC^\times$-orbits, 
\item[hyp] if $1<\kappa_0<m_\hyp$ for some $m_\hyp>1$, then $h$ is nondegenerate hyperbolic and such that the tangent spaces to the $\CC^\times$-orbits are negative.
\end{description}
This implies that  $\PP(V^\circ)$ acquires a geometric structure which is respectively modeled
on Fubini-Study space, flat complex Euclidean space and complex hyperbolic space. A suitable combination of rationality and symmetry conditions (which generalizes the half integrality condition \ref{def:half}) yields a generalization of Theorem \ref{thm:main}. We thus obtain new examples of groups operating discretely and with finite covolume on a complex ball.

\end{document}